\documentclass[twoside]{article}
\usepackage{amsmath, amsthm, amsfonts, amssymb, latexsym}
\usepackage{mathptmx}
\usepackage{graphicx}
\usepackage{fancybox}
\usepackage{color}

\numberwithin{equation}{section}

\Large

\def\a{\alpha}
\def\e{\varepsilon}

\def\R{\bf R}

\newtheorem{thm}{Theorem}[section]
\newtheorem{df}{Definition}[section]

\newtheorem{lem}{Lemma}[section]
\newtheorem{cor}{Corollary}[section]
\newtheorem{prop}{Proposition}[section]
\newtheorem{rem}{Remark}[section]
\numberwithin{equation}{section}
\newtheorem{exam}{Example}[section]
\begin{document}

\title{ Hardy's inequalities  \\ with non-doubling weights and   sharp remainders. }
\author{    Toshio Horiuchi   }
\date{}
\maketitle

\begin{abstract}
In the present paper 
we shall establish     $N$-dimensional   Hardy's inequalities with non-doubling
  weight functions of the  distance $\delta(x)$ to the   boundary $\partial\Omega$, where 
  $\Omega$   is a $C^2$ class bounded domain of $\mathbf R^N\, ( N\ge 1)$.
This work is essentially based on one  dimensional  weighted Hardy's inequalities with  one-sided boundary condition and 
 sharp remainders. 
As weights we  admit  rather general ones that  may vanish or blow up in infinite order
 such as  $e^{-1/t}$  or  $e^{1/t}$ at $t=0$ in one-dimensional case.
  \footnote{ 2010 mathematics Subject Classification
{Primary 35J70, Secondary 35J60, 34L30, 26D10 }; keywords:{ Weighted Hardy's inequalities, Weak  Hardy property, $p$-Laplace operator with weights}; This research was partially supported by Grant-in-Aid for Scientific Research (No. 20K03670, No. 21K03304).
}
 \end{abstract}

\section{ Introduction }
In the  present paper, 
we shall  begin with introducing one-dimensional  weighted Hardy's inequalities with sharp remainders under one-sided boundary condition. 
As weights  we shall deal with the  so-called   non-doubling  weights in addition to usual doubling ones.
 Then we shall establish  $N$-dimensional weighted Hardy's inequalities with non-doubling
  weights of the  distance  $\delta(x):= { \rm dist}(x,\partial\Omega)$ to the   boundary $\partial\Omega$, where 
  $\Omega$   is a $C^2$ class bounded domain of $\mathbf R^N\, ( N\ge 2)$. 
  A  positive \textcolor{black}{ continuous} function $w(t)$ on $(0,\infty)$ is  said to be a  doubling weight  if  there exists a positive number $C$  such  that  we  have 
\begin{equation} C^{-1} w(t)\le w(2t)\le C w(t) \quad (0<t<\infty),\label{doubling} \end{equation}
where $C$ is independent of  each $t \in (0,\infty)$.
When $w(t)$ does not possess this property, $w(t)$ is  said to be  a non-doubling weight in the  present paper.
In one-dimensional case we   typically treat a weight function $w(t)$ that  may vanish or  blow up in infinite order such  as $e^{-1/t}$  or  $e^{1/t}$ at $t=0$. In  such cases the limit of ratio $w(t)/w(2t)$ as  $t\to +0$ may become $0$ or $+\infty$, and hence they are regarded as non-doubling weights according to our notion.
\begin{df} Let $1<p<\infty$  and 
set $p'= p/(p-1)$.  $\Lambda_p$  denotes the  one-dimensional Hardy best constant defined by
$$\Lambda_p=\frac{1}{(p')^p}=\left(1- \frac1p\right)^{p}.$$

\end{df}
Let   $\eta>0$. By  $ C_c^1((0,\eta])$ we denote the  set of  all $C^1$  functions  with  compact  supports in $(0,\eta]$. 
Then
one-dimensional Hardy's inequalities with  one-sided boundary condition in this  paper  are typically represented by the followings.
\begin{prop}\label{Prop} Assume that  $1<p<\infty$,   $\mu>0$ and $\eta>0$.  \par\noindent
1.
 For  every  $u\in C_c^1((0,\eta])$ we  have
\begin{equation} \int_0^\eta |u'(t)|^p e^{-(p-1)/t}\,dt + \frac{(\Lambda_{p})^{1/p'} }{\mu^{p-1} }|u(\eta)|^p\ge
\Lambda_{p} \int_0^\eta \frac{|u(t)|^p {e^{-(p-1)/t}}\,dt}{\left(e^{-1/t}\left( \int_t^\eta e^{1/s}\,ds +\mu \right)\right)^p}. \label{DD}
\end{equation}
\par \noindent 2.
For  every  $u\in C_c^1((0,\eta])$ we  have
\begin{equation} \int_0^\eta |u'(t)|^p e^{(p-1)/t}\,dt \ge 
\Lambda_{p} \int_0^\eta \frac{|u(t)|^p {e^{(p-1)/t}}\,dt}{\left( e^{1/t} \int_0^t e^{-1/s} \,ds \right)^p}+ \frac{(\Lambda_{p})^{1/p'}}{\left(\int_0^\eta e^{-1/s}\,ds\right)^{p-1} }|u(\eta)|^p
. \label{GG}
\end{equation}

\end{prop}
Proposition \ref{Prop}   will be  established  as Corollary \ref{cor2.2} to Theorem \ref{T3.1}.
When   $u(\eta)=0$,  these inequalities  are a variant of classical weighted Hardy's inequalities (see \cite{ KO,ahl,Ma}).
It is interesting that both  coefficients $\Lambda_{p}$ and $(\Lambda_{p})^{1/p'}$ appearing in (\ref{DD}) and (\ref{GG}) are best, even though the  inequalities contain   two parameters  $\eta$ and $\mu$.
To  see    the  sharpness of (\ref{DD}) and (\ref{GG}), by  the  density argument it  suffices  to employ $u_\varepsilon(t)= \left( \mu+ \int_t^\eta e^{1/s}\,ds\right)^{1/p'-\varepsilon}$ and 
$u_\varepsilon(t)= \left(  \int_0^t e^{-1/s}\,ds\right)^{1/p'+\varepsilon}$ for
   test  functions ($\varepsilon \to +0$) respectively. (For the detail see Subsection 4.2 (Part 1).)
\par
Our  first  purpose  in  this  paper is not  only  to establish a general version of Proposition \ref{Prop}  but  also  improve  it by  adding sharp remainder terms.
By $W(\mathbf R_+)$  we  denote  a class of  functions 
$$ \{ w\in C^1({\R}_+): w>0, \lim_{t\to+0}w(t)=a\, \text{  for  some }\,a\in [0,\infty] \}$$  
with $\mathbf R_+=(0,\infty)$.
As weights   we adopt  functions $W_p(t)=w(t)^{p-1}$ with $w(t)\in P(\mathbf R_+) \cup Q(\mathbf R_+)$, where   
 \begin{equation}  \begin{cases}&P(\mathbf R_+)= \{ w(t)\in W(\mathbf R_+) : \,  w(t)^{-1} \notin L^1((0,\eta)) \, \text{ for some} \, \eta >0\},\\
 &Q(\mathbf R_+) =\{ w(t)\in W(\mathbf R_+) :  \, w(t)^{-1}\in L^1((0,\eta)) \, \text{ for any } \, \eta >0 \}.
 \end{cases}
 \end{equation}
Clearly   $W(\mathbf R_+)=P(\mathbf R_+) \cup Q(\mathbf R_+)$,
$e^{-1/t} \in  P(\mathbf R_+)$  and  $e^{1/t} \in Q(\mathbf R_+)$ for $t>0$ ( For the precise definitions see Section 2 ).
By  virtue  of Proposition \ref{Prop}, it  is  clearly  seen  that
our   results on  this  matter   essentially depend on whether $w$ belongs to
 $ P(\mathbf R_+) $ or $ Q(\mathbf R_+)$.
In particular when $w(t)\in P(\mathbf R_+)$, it  follows from Proposition \ref{ct1} that
\begin{equation}
\inf_{u\in V_\eta}
\int_0^\eta  |u'(t)|^p W_p(t)\,dt =0, \label{C}
 \end{equation}
 where $V_\eta= \{  u\in C^1([0,\eta]) :   u(0)=0, u(\eta)=1\}.  $
Nevertheless   we   have   sharp Hardy type  inequalities  (\ref{nct1}) and  (\ref{nct2}) in Theorem \ref{T3.1}.
\par
As an important application, in Theorem \ref{CT3}  we  shall establish  
 $N$-dimensional  Hardy's inequalities with 
  weights  being  functions  of the  distance  $\delta (x)=\rm{ dist}(x,\partial\Omega)$ to the   boundary $\partial\Omega$. 
   In  this  task it  is  crucial to establish sharp weighted  Hardy's  inequalities in  the tubler neighborhood $\Omega_\eta$  of $\Omega$, which are reduced  to  the  one  dimensional inequalities  in Theorem \ref{T3.1}.
To  this  end,   $\Omega $ is  assumed  to be  a bounded domain of  $\mathbf R^N$ ( $N\ge 2$ ) whose boundary  $\partial\Omega$ is a $C^2$ compact manifolds in  the  present  paper. 
We  prepare  more  notations to  describe  our  results.
For  $ W_p(t) = w(t)^{p-1}$ with $ w(t)\in W(\mathbf R_+)$, we  define  a weight function $W_p(\delta(x))$ on $\Omega$  by $$W_p(\delta(x))=(W_p\circ\delta)(x).$$
By
  $ L^p(\Omega; W_p(\delta))$ we denote the space of Lebesgue measurable functions with weight $ W_p(\delta(x))$, 
for which 
\begin{equation} \| u \|_{ L^p(\Omega; W_p(\delta))} = \bigg( \int_{\Omega}|u(x)|^p W_p(\delta(x))\,dx\bigg ) ^{1/p} < +\infty.\label{2.1} \end{equation} 
By  $ C_c^\infty(\Omega)$ we denote the  set of  all $C^\infty$  functions  with  compact  supports in $\Omega$. 
$W_{0}^{1,p}(\Omega,W_p(\delta))$ is given by the  completion of $C_c^\infty(\Omega)$ with respect to the norm 
defined by
\begin{equation} 
\| u\|_{W^{1,p}_{0}(\Omega; W_p(\delta)) } =
\| |\nabla u|  \|_{ L^p(\Omega; W_p(\delta))} + \| u\|_{ L^p(\Omega; W_p(\delta))}. \label{2.2}
\end{equation}
Then 
 $W^{1,p}_{0}(\Omega; W_p(\delta)) $  becomes a Banach space  with the norm $\| \cdot\|_{W^{1,p}_{0}(\Omega; W_p(\delta)) }
$.
 Under these preparation 
we will  establish $N$-dimensional weighted Hardy's inequality as
Theorem \ref{CT3}, which is    the   counter-part  to  Theorem \ref{T3.1}.
In particular for $w(t)\in Q(\mathbf R_+)$,
as in Corollary  \ref{NCC1}  we   have a simple inequality which  is a generalization of classical Hardy's inequality:
\begin{equation}\label{HI}
\int_\Omega  |\nabla u|^p W_p(\delta(x)) \,dx\ge \gamma \int_\Omega \frac{|u|^pW_p(\delta(x))}{F_\eta(\delta(x))^p}\,dx, \quad \forall u\in W^{1,p}_{0}(\Omega; W_p(\delta)) , 
 \end{equation}
where $\eta$ is a sufficiently small positive number, $\gamma$  is some  positive  constant and $F_\eta(\delta(x))=(F_\eta\circ \delta)(x)$ is a nonnegative function defined in Definition \ref{D2}. 
If $w=1$, then $F_{\eta}(t)=t \, (t\le \eta)$; $\eta\,(t\ge \eta)$ and 
  (\ref{HI})  is  a well-known Hardy's inequality \textcolor{black}{having $\left(\min(\delta(x),\eta)\right)^{-p}$ as the  Hardy potential,} which  is valid for  a bounded domain $\Omega $  of  $\mathbf R^N$ with Lipschitz boundary  (cf. \cite{BM,D, MMP}).
Further
if     $\Omega$ is  convex,  then   $\gamma = \Lambda_{p}$ holds  for   arbitrary $1<p<\infty$ (cf. \cite{MS,YZ}). \par 
It  is  worthy to remark that (\ref{HI}) is  never valid in  the    case that $w(t)\in P(\mathbf R_+)$
by (\ref{C}) (see also  Proposition \ref{ct1} and Proposition \ref{CT1}). Nevertheless,   in  this   case
 we  shall   establish   weighted  Hardy's inequalities with a switching function in  Theorem \ref{CT3} and  Corollary \ref{NCC1}, which correspond to  Theorem \ref{T3.1} and  its corollaries. 
We remark  that these   Hardy's  inequalities  with  a compact perturbation  are  closely relating   to  
 the so-called weak Hardy property of  $\Omega$.  In  fact, if $w(t) \in P(\mathbf R_+)$, then 
 a constant $\gamma^{-1}$ in (\ref{2.7})  concerns  the weak Hardy constant,  but  in  this  case   the  strong Hardy  constant  is $+\infty$ (  see \cite{D} for  the  detail). In \cite{ahl},  we  have improved     Hardy's inequalities adopting  $\delta(x)^{\a p}$  ( powers  of  the distance $\delta(x)$ to the boundary $\partial\Omega$ )   as   weight  functions   instead of  $W_p(\delta)$.  In  the  present paper,  some  inequalities 
of  Hardy type  in \cite{ah0} and \cite{ahl} will  be   employed with minor  modifications, especially  when $1<p<2$
 (see also \cite{anm}).\par
We remark  that our  results   will  be applicable  to variational problems with critical Hardy potentials in a  coming paper \cite{ah4} (c.f. \cite{ah3}) 
and  also applicable to the Caffarelli-Kohn-Nirenberg type inequalities with non-doubling weights in the coming paper \cite{ah5}. 
 \par
 This  paper is  organized in  the  following way: In Section 2 we introduce a  class of weight functions  $W(\mathbf R_+)$ and two subclasses $P(\mathbf R_+)$  and $ Q(\mathbf R_+)$ together with 
 functions such as $F_\eta(t)$ and $G_\eta(t)$, 
which are crucial  in this  paper. Further a notion of  admissibilities for $P(\mathbf R_+)$  and $ Q(\mathbf R_+)$ is  introduced.
In Section 3, the main  results are described.  
The results are divided into two cases (one-dimensional case and $N$-dimensional case), which are described in Subsection 3.1 and Subsection 3.2, respectively.
Theorem \ref{T3.1} and Theorem \ref{T3.2} are  established in Section 4.  Theorem \ref{CT3}  together with  Corollary \ref{NCC1}  are proved in Section 5. 
The  proof  of Theorem \ref{NCT2} is given in Section 6 and the  proofs of  Proposition \ref{ct1}  and Proposition \ref{CT1} are  given  in Section 7.
In  Appendix  the  proof of Lemma \ref{lemma4.4} is provided.
 Some auxiliary inequalities are also given as Lemma \ref{l2}.
 \par\medskip
 \section{Preliminaries  }

First we  introduce a class of weight functions which  is crucial in this paper.
\begin{df}\label{D0}
 Let us  set $\mathbf R_+= (0,\infty)$ and 
 \begin{equation}W(\mathbf R_+)=
 \{ w(t)\in C^1({\R}_+): w(t)>0, \lim_{t\to+0}w(t)=a\, \text{  for  some }\,  a\in [0,\infty] \}.
\end{equation}

\end{df}
In  the  next  we define two subclasses of  this  rather large space.

\begin{df}\label{D1} Let  us  set 
 \begin{equation} P(\mathbf R_+)= \{ w(t)\in W(\mathbf R_+) : \,  w(t)^{-1} \notin L^1((0,\eta)) \, \text{ for some} \, \eta >0\}.\end{equation}
 \begin{equation} 
 Q(\mathbf R_+) =\{ w(t)\in W(\mathbf R_+) :  \, w(t)^{-1}\in L^1((0,\eta)) \, \text{ for any } \, \eta >0 \}.\end{equation}
\end{df}

\begin{rem} \begin{enumerate}
\item From  Definition \ref{D0} and Definition \ref{D1} it  follows   that $W(\mathbf R_+)= P(\mathbf R_+)\cup Q(\mathbf R_+)$ and  $P(\mathbf R_+)\cap Q(\mathbf R_+)=\phi.$
\item If  $w^{-1} \notin L^1((0,\eta))  \, \text{ for some} \, \eta >0$, then  $w^{-1} \notin L^1((0,\eta)) \, \text{ for any} \, \eta >0$. Similarly
if $w^{-1} \in L^1((0,\eta))  \, \text{ for some} \, \eta >0$, then  $w^{-1} \in L^1((0,\eta)) \, \text{ for any} \, \eta >0$.
\item
If  $w\in  P(\mathbf R_+)$, then $\lim_{t\to +0}w(t)=0$. Hence 
by  setting $w(0)=0$, $w$ is uniquely  extended to  a continuous function on $[0,\infty)$. 
On the other hand  if $w\in  Q(\mathbf R_+)$, then    possibly  $\lim_{t\to +0}w(t)=+\infty$.
\end{enumerate}
\end{rem}
Here we give some  fundamental examples:
\begin{exam} 
\begin{enumerate}\item $ e^{-1/t}\in P(\mathbf R_+)$ and $ e^{1/t}\in Q(\mathbf R_+)$. (See Corollary \ref{cor2.2}.)
\item For $p'=p/(p-1)$,
$t^{\a p'} \in P(\mathbf R_+)$ if $\a \ge 1/p'$  and 
$t^{\a p'} \in Q(\mathbf R_+)$ if $\a < 1/p'$. (See Corollary \ref{cor2.3} with $W_p(t)=w(t)^{p-1}= t^{\a p}$.)
\item For $ \a\in \mathbf R$, $ t^\a e^{-1/t}\in P(\mathbf R_+)$ and $t^\a e^{1/t}\in Q(\mathbf R_+)$.
\end{enumerate}
\end{exam}

In  the next we define functions such as $F_\eta(t)$ and $G_\eta(t)$, 
which are crucial in considering variants of the Hardy potential like $F_\eta(\delta(x))^{-p}$ in (\ref{HI}).
\begin{df}\label{D2} Let $\mu>0$ and   $\eta>0$. For
 $w \in W(\mathbf R_+)$, we  define the followings:
 \begin{enumerate}
 \item  When $w\in P(\mathbf R_+)$, 
   \begin{align}
 F_{\eta}(t; w,\mu) & = \begin{cases} &w(t)\left( \mu+ \int_t^\eta {w(s)}^{-1}\, {ds}\right), \qquad   \text{ if } \, t\in (0,\eta], \\ 
 &w(\eta) \mu, \qquad \qquad \qquad \qquad\quad    \text{ if } \, t\ge  \eta.
 \end{cases}\\
 G_\eta(t;w,\mu) &= \begin{cases} &\mu+  \int_t^\eta  {F_\eta(s;w,\mu)}^{-1}\, {ds}, \qquad     \text{ if } \, t\in (0,\eta],\\
&\mu, \qquad \qquad \qquad \qquad\qquad \quad \text{ if } \, t\ge \eta.
 \end{cases} \label{PG}
\end{align}

\item When $w\in Q(\mathbf R_+)$, 
 \begin{align}
 F_{\eta}(t; w) &= \begin{cases} &w(t)\int_0^t {w(s)}^{-1} \,{ds}, \qquad \qquad   \text{ if } \, t\in (0,\eta], \\ 
 &w(\eta) \int_0^\eta w(s)^{-1} \,ds, \qquad \quad           \text{ if } \, t\ge \eta.
\end{cases}\\
  G_\eta(t;w,\mu)&= \begin{cases} &\mu+  \int_t^\eta  {F_\eta(s;w)}^{-1}\, {ds}, \qquad    \text{ if } \, t\in (0,\eta], \\
&\mu, \qquad \qquad \qquad \qquad \qquad  \text{ if } \, t\ge \eta.
 \end{cases}
\end{align}
\item $F_{\eta}(t;w,\mu)$ and $F_{\eta}(t;w)$ are abbreviated as $F_{\eta}(t)$.
$ G_{\eta}(t; w,\mu)$ is abbreviated as  $ G_{\eta}(t)$.
 \item  For   $w\in W(\mathbf R_+)$,
 we  define
 \begin{equation}  W_p(t)= w(t)^{p-1}. \label{1.4}\end{equation}
\end{enumerate}
\end{df}

\begin{rem}\label{defofG} In the definition (\ref{PG}), one can replace $G_\eta(t; w,\mu)$ with the more  general 
$G_\eta(t; w,\mu,\mu')$ $= \mu'+  \int_t^\eta  {F_\eta(s;w,\mu)}^{-1}\, {ds},  \, t\in (0,\eta]; {\, } \mu', \, t\ge \eta \,$ 
for $\mu'>0 $.
However, for simplicity  this  paper uses  (\ref{PG}).  
\end{rem}

\begin{exam} Let $w(t)=t^{\a p'} $ for $\a\in \R$, $1<p<\infty$ and $p'=p/(p-1)$. 
\begin{enumerate}
\item When $\a>1/p'$, $F_\eta(t)=t/(\a p'-1) \,$   and 
$G_\eta(t)=\mu +(\a p'-1)\log(\eta/t)$ provided that  $ \mu= \eta^{1-\a p'}/(\a p'-1)$.
\item When $\a=1/p'$,  $F_\eta(t)=t(\mu+ \log(\eta/t)) \,$   and 
$G_\eta(t)=\mu -\log \mu +\log\left( \mu+ \log(\eta/t)\right)$.
\item When $\a<1/p'$, $F_\eta(t)=t/(1-\a p') \,$   and 
$G_\eta(t)=\mu +(1-\a p')\log(\eta/t)$.
\end{enumerate}
\end{exam}
By using integration by parts  we  see  the followings:
\begin{exam}
\begin{enumerate}\item If either $w(t)=e^{-1/t} \in P(\R_+)$ or $w(t)=e^{1/t}\in Q(\R_+)$, then $F_\eta(t)= O(t^2)$ as $t\to +0$.
\item Moreover,
if $w(t)= \exp({\pm t^{-\a}}) $ with $\a> 0$, then $F_\eta(t)= O(t^{\a+1})$ as $t\to +0$. 
In fact, it  holds  that $\lim_{t\to +0} F_\eta(t)/t^{\a+1}=1/\a$. 
\end{enumerate}
\end{exam}

In a similar way  we define the following:
\begin{df}\label{fg} Let $p'= p/(p-1)$, $\mu>0$ and   $\eta>0$. For
 $w \in W(\mathbf R_+)$  and $t\in (0,\eta)$, we  define the followings:
 \begin{enumerate}
 \item When $w\in P(\mathbf R_+)$,
  \begin{align}
f_{\eta}(t; w,\mu) &= \begin{cases} & \mu+ \int_t^\eta {w(s)}^{-1}\, {ds}, \qquad   \text{ if } \, t\in (0,\eta],\\ 
 &\mu, \qquad \qquad     \qquad \qquad       \text{ if } \, t\ge \eta.
\end{cases}\\
g_\eta(t;w,\mu)&=  (p' f_\eta(t;w,\mu))^{{1}/{p'}}, \qquad \qquad           \text{ if } \, t>0.
\end{align}
 \item When $w\in Q(\mathbf R_+)$,
  \begin{align}
f_{\eta}(t; w) &= \begin{cases}&
\int_0^t {w(s)}^{-1} \,{ds}, \qquad \qquad           \text{ if } \, t\in (0,\eta],\\
&\int_0^\eta w(s)^{-1}\,ds,\qquad \qquad           \text{ if } \, t\ge \eta.
\end{cases}\\
 g_\eta(t;w) &= (p' f_\eta(t;w))^{{1}/{p'}}, \qquad \qquad  \quad         \text{ if } \, t>0.
\end{align}
\item $f_{\eta}(t;w,\mu)$ and $f_{\eta}(t;w)$ are abbreviated as $f_{\eta}(t)$.
$g_{\eta}(t;w,\mu)$ and $g_{\eta}(t;w)$ are abbreviated as  $ g_{\eta}(t)$.

\end{enumerate}
\end{df}

\begin{rem}\label{Remark2.2}\begin{enumerate}\item
We note that for  $0<t<\eta$
\begin{equation}\label{2.13}\begin{cases}& \frac {d}{dt}\log f_{\eta}(t)= -{F_\eta(t) }^{-1}, \quad \text{ if } w\in P({\mathbf R}_+),\\ 
&\frac {d}{dt}\log f_{\eta}(t)= {F_\eta(t) }^{-1},\quad  \mbox{ if } w\in Q({\mathbf R}_+),\\
&\frac {d}{dt}\log G_{\eta}(t)= -({F_\eta(t) G_{\eta}(t)})^{-1},\quad  \mbox{ if } w\in W({\mathbf R}_+),\\
&\frac {d}{dt}G_{\eta}(t)^{-1}= ({F_\eta(t) G_{\eta}(t)^2})^{-1},\quad  \mbox{ if } w\in W({\mathbf R}_+).
\end{cases}
\end{equation}
By Definition \ref{D1} and (\ref{2.13}), we  see that $F_\eta(t)^{-1}\notin L^1((0,\eta))$, $( F_\eta (t)G_\eta(t))^{-1}\notin L^1((0,\eta))$ but $ (F_\eta(t) G_\eta(t)^2 )^{-1}\in L^1((0,\eta))$.   
\item
If $w\in W(\R_+)$, then immediately  we  have $\liminf_{t\to +0}F_\eta(t)=0$ from 1. 
\end{enumerate}
\end{rem}
\begin{exam}
If either
$w(t)=t^2e^{-1/t}\in P(\mathbf R_+)$ or $w(t)=t^2e^{1/t}\in Q(\mathbf R_+)$, then
 $F_{\eta}(t)=O(t^2)$ and  $G_{\eta}(t)=O(1/t)$ as $t\to 0$.
\end{exam}
Now  we introduce  two admissibilities for $P(\mathbf R_+)$ and  $Q(\mathbf R_+)$.
\begin{df} \begin{enumerate}\item
A function $w\in P(\mathbf R_+)$ is said  to be admissible if there exist  positive  numbers $\eta$  and $K$ such  that we  have 
 \begin{equation}\label{cond1}
 \int_t^\eta w(s)^{-1} \,ds \le e^{K/\sqrt t} \qquad (0<t<\eta).
 \end{equation}
 \item 
 A function $w\in Q(\mathbf R_+)$ is said  to be admissible if there exist  positive  numbers $\eta$  and $K$ such  that we  have 
\begin{equation}\label{cond2}
 \int_0^t w(s)^{-1} \,ds \ge e^{-K/\sqrt t} \qquad (0<t<\eta).
 \end{equation}
 \end{enumerate}
  \end{df}
  
\begin{df}  \label{admissibility}
By $P_A (\mathbf R_+)$  and $ Q_A(\mathbf R_+)$ we  denote  the  set of  all  admissible  functions  in $P(\mathbf R_+)$ and $ Q(\mathbf R_+)$ respectively. We  set \begin{equation}W_A(\mathbf R_+)= P_A (\mathbf R_+)\cup Q_A(\mathbf R_+).\end{equation}
\end{df}
Here we give  typical examples:
\begin{exam} $e^{-1/t}\notin P_A(\mathbf  R_+)$, $e^{1/t}\notin Q_A(\mathbf  R_+)$ but $e^{-1/\sqrt{t}}\in P_A(\mathbf  R_+)$, $e^{1/\sqrt{t}}\in Q_A(\mathbf  R_+)$.\par\medskip
\par\noindent
{\bf Verifications: } \par
\noindent  $e^{-1/t} \notin  P_A(\R_+):$ For  small $t>0$, we  have $\int_t^\eta e^{1/s} \,ds \ge \int_t^{2t} e^{1/s} \,ds\ge te^{1/(2t)}$.
But this contradicts to (\ref{cond1}) for  any $K>0$. \par\noindent
{ $e^{-1/\sqrt t} \in  P_A(\R_+) :$}
Since $e^{1/\sqrt s}\le e^{1/\sqrt t}\, (t<s<\eta)$, we  have $\int_t^\eta e^{1/\sqrt s}\,ds \le \eta  e^{1/{\sqrt t}}\le e^{K/\sqrt t}$ for  some $K>1$. 
\par\noindent $e^{-1/t} \notin  Q_A(\R_+):$
For $0<s\le t$, we  have $\int_0^t e^{-1/s}\,ds\le t e^{-1/t}$. But this contradicts to (\ref{cond2}) for  any $K>0$.
 \par\noindent
{$e^{-1/\sqrt t} \in Q_A(\R_+):$} For $t/2<s<t$, we  have 
$\int_0^t e^{-1/\sqrt s}\,ds \ge \int_{t/2}^t e^{-1/\sqrt s}\,ds$ $\ge (t/2) e^{-\sqrt 2/ \sqrt s}$ $\ge e^{-K/\sqrt t}$ for some $ K>\sqrt 2$.

\end{exam}
\begin{prop}\label{prop2.1}
\begin{enumerate}\item
Assume that either   $w\in P_A (\mathbf R_+)$ or  $w\in Q_A (\mathbf R_+)$.  Then 
there exist positive numbers $\eta$ and $K$ such  that we  have 
\begin{equation}
\sqrt t \,G_\eta (t) \le K\qquad  t\in (0,\eta).\label{2.10}
\end{equation}
\item 
Assume that    $w\in W(\mathbf R_+)$  and $w$ satisfies  (\ref{2.10}) for  some positive numbers $\eta$ and $K$. 
Then 
$w$ is  admissible in  the sense of Definition \ref{admissibility}, that is,  $w$ belongs  to $W_A (\mathbf R_+)$.
\end{enumerate}
\end{prop}

\noindent {\bf Proof:} 

By integrating
(\ref{2.13}), we have   for  $0<t<\eta$
\begin{equation}\begin{cases}&G_\eta (t)= \log \left( \int_t^\eta { w(s)}^{-1}\,ds +\mu\right)-\log \mu + \mu , \quad \text{ if } w\in P(\mathbf R_+),\\ \\
&G_\eta (t)= \log \left( \int_0^\eta { w(s)}^{-1}\,ds \right)- \log\left(\int_0^t { w(s)}^{-1}\,ds \right)  + \mu ,\quad  \mbox{ if } w\in Q(\mathbf R_+).
\end{cases}
\end{equation}
Hence  the  inequality  (\ref{2.10})  with   positive  numbers $\eta$  and $K$ is  equivalent  to 
\begin{equation}\begin{cases}& \int_t^\eta { w(s)}^{-1}\,ds\le \mu( e^{ {K}/{\sqrt t} -\mu}-1)\quad (0<t<\eta), \quad \text{ if } w\in P_A(\mathbf R_+),\\ \\
&\int_0^t {w(s)}^{-1}\,ds \ge e^\mu \int_0^\mu  {w(s)}^{-1}\,ds \,e^{-{K}/{\sqrt t}}\quad    (0<t<\eta),\quad  \mbox{ if } w\in Q_A(\mathbf R_+).
\end{cases} \label{2.10'}
\end{equation}
Here we  
note  that  for each $\mu>0$, $\eta>0$ and  $w\in Q(\mathbf R_+)$  there exist  some positive numbers  $K_1$, $K_2$ and $K_3$  such  that  $K_1\le K_2\le K_3 $ and we  have

\begin{equation} \begin{cases}&
e^{ {K_1}/{\sqrt t} } \le  \mu( e^{ {K_2}/{\sqrt t} -\mu}-1) \le e^{ {K_3}/{\sqrt t}} \quad (0<t<\eta)\\ &\\
&
e^{- {K_3}/{\sqrt t}}\le e^\mu \int_0^\mu  {w(s)}^{-1}\,ds \,\,e^{-{K_2}/{\sqrt t}} \le 
e^{- {K_1}/{\sqrt t}} \quad    (0<t<\eta).\end{cases}
\end{equation}
\noindent Therefore the assertion is now clear. \hfill$ \Box$
\section{Main results} 
\subsection{Results in the  one  dimensional  case}
We introduce  function spaces  to state the results in one dimensional case:\par\noindent
By
  $ L^p((0,\eta]; W_p)$ we denote the space of Lebesgue measurable functions with weight $ W_p(t)$, 
for which 
\begin{equation} \| u \|_{ L^p((0,\eta]; W_p)} = \bigg( \int_0^\eta|u(t)|^p W_p(t)\,dt\bigg ) ^{1/p} < +\infty.\label{2.1'} \end{equation} 
By  $ C_c^\infty((0,\eta])$ we denote the  set of  all $C^\infty$  functions  with  compact  supports in $(0,\eta]$. 
 $W_{0}^{1,p}((0,\eta]; W_p)$ is given by the  completion of $C_c^\infty((0,\eta])$ with respect to the norm 
defined by
\begin{equation} 
\| u\|_{ W_{0}^{1,p}((0,\eta]; W_p)} =
\|u' \|_{ L^p((0,\eta]; W_p)} + \| u\|_{ L^p((0,\eta]; W_p)}. \label{2.2'}
\end{equation}
Then 
 $W_{0}^{1,p}((0,\eta]; W_p) $  becomes a Banach space  with the norm $\| \cdot\|_{ W_{0}^{1,p}((0,\eta]; W_p)}$. 
 We also define a switching function:
\begin{df} $($ \rm{Switching function} $)$ 
For $w\in W(\mathbf R_+)= P(\mathbf R_+) \cup Q(\mathbf R_+)$ we  set
\begin{equation} s(w)= \begin{cases} -1 \quad & \text{if } \quad w\in P(\mathbf R_+),\\
\,\, 1\quad & \text{if } \quad w\in Q(\mathbf R_+).\end{cases}
\end{equation}
\end{df}
We state one-dimensional Hardy's inequalities with  one-sided boundary condition. 
 \par\medskip
\begin{thm}\mbox{$(   w\in W(\mathbf R_+) )$ } \label{T3.1}\par \noindent
Assume that $1<p<\infty$, $\eta>0$, $\mu>0$ and  $w\in W(\mathbf R_+)= P(\mathbf R_+) \cup Q(\mathbf R_+)$. Then we have  the  followings:
\begin{enumerate}\item For  every $u\in  W_{0}^{1,p}((0,\eta]; W_p)\cap C((0,\eta])$, we have
\begin{align}  \int_0^\eta   |u'(t)|^p W_p(t)\,dt \ge \Lambda_{p} 
 \int_0^\eta \frac{ |{u(t)}|^p W_p(t)\,dt }{F_\eta(t)^p}+s(w)\frac{(\Lambda_p)^{1/p'}}{f_\eta(\eta)^{p-1}} |u(\eta)|^p
.\label{nct1}
\end{align}
Moreover the coefficients $\Lambda_p$  and $s(w)(\Lambda_p)^{1/p'}$ are optimal.
\item There exist  positive numbers  $C=C(w,p,\eta,\mu )$ and  $L=L(w,p,\eta,\mu )$ such  that  for  every 
$u\in W_{0}^{1,p}((0,\eta]; W_p)\cap C((0,\eta])$, we  have 
\begin{equation}
\begin{split}  \int_0^\eta   |u'(t)|^p W_p(t)\,dt  
\ge & \Lambda_{p} 
 \int_0^\eta \frac{ |{u(t)}|^p W_p(t)\,dt }{F_\eta(t)^p} + C\int_0^\eta \frac{ |u(t)|^p W_p(t) \,dt}{ F_\eta(t)^p G_\eta(t)^2}
 +s(w)L|u(\eta)|^p,\label{nct2}
\end{split}
\end{equation}
where $C$ and $L$ can be taken  independent of each $u$.
\end{enumerate}
\end{thm}
The following is a direct consequence from this theorem.
We  remark that $C_c^1((0,\eta])$ is  densely contained in $W_{0}^{1,p}((0,\eta]; W_p) $.
\begin{cor}\label{cor2.2}
\begin{enumerate}\item Let  $1<p<\infty$, $\eta>0$  and  $\mu>0$.  
Then,  for  every  $u\in C_c^1((0,\eta])$ we  have
\begin{equation} \int_0^\eta |u'(t)|^p e^{-(p-1)/t}\,dt + \frac{(\Lambda_{p})^{1/p'}}{\mu^{p-1} }|u(\eta)|^p\ge
\Lambda_{p} \int_0^\eta \frac{|u(t)|^p {e^{-(p-1)/t}}\,dt}{\left(e^{-1/t}\left( \int_t^\eta e^{1/s}\,ds +\mu \right)\right)^p}. \label{D}
\end{equation}
\item
Let  $1<p<\infty$  and  $\eta>0$.
Then,   for  every  $u\in C_c^1((0,\eta])$ we  have
\begin{equation} \int_0^\eta |u'(t)|^p e^{(p-1)/t}\,dt \ge 
\Lambda_{p} \int_0^\eta \frac{|u(t)|^p {e^{(p-1)/t}}\,dt}{\left( e^{1/t} \int_0^t e^{-1/s} \,ds \right)^p}+ \frac{(\Lambda_{p})^{1/p'}}{\left(\int_0^\eta e^{-1/s}\,ds\right)^{p-1} }|u(\eta)|^p
. \label{G}
\end{equation}
\end{enumerate}
\end{cor}
\begin{df}For $1<p<+\infty$ and $\a \in \mathbf R$  we set
\begin{equation}\Lambda_{\a,p}=\left |\frac{1}{p'}-\a\right |^p = \left |\frac{p-1-\a p}{p}\right |^p.\label{2.5}
\end{equation} 
\end{df}

\begin{cor}\label{cor2.3}
 \begin{enumerate}\item If $\a>1/p'$, $1<p<\infty$ and $\eta>0$,  then for  every  $u\in C_c^1((0,\eta])$ we  have
\begin{equation} \int_0^\eta |u'(t)|^p t^{\a p}\,dt + \frac{(\Lambda_{\a,p})^{1/p'}}{\eta^{p-1-\a p}} |u(\eta)|^p\ge\Lambda_{\a,p} \int_0^\eta {|u(t)|^p}{t^{(\a -1)p}}\,dt. \label{F}
\end{equation}
 \item 
If $\a=1/p'$, $1<p<\infty$, $\eta>0$ and $R>e$,  then  for  every  $u\in C_c^1((0,\eta])$ we have
\begin{equation} \int_0^\eta |u'(t)|^p t^{p-1}\,dt 
+ \frac{(\Lambda_{p})^{1/p'} } { (\log R)^{p-1} } |u(\eta)|^p \ge\Lambda_{p} \int_0^\eta \frac{|u(t)|^p}{t A_1(t/\eta)^p}\,dt,\label{E}
\end{equation} 
where $ A_1(t)= \log (R/t)$ with $R>e$.
\item
If $ \a<1/p'$, $1<p<\infty$ and $\eta>0$, then for  every  $u\in C_c^1((0,\eta])$ we  have
\begin{equation} \int_0^\eta |u'(t)|^p t^{\a p}\,dt \ge 
\Lambda_{\a,p}\int_0^\eta {|u(t)|^p}{t^{(\a-1) p}}\,dt +
\frac{({\Lambda_{\a,p}})^{1/p'}}{\eta^{p-1-\a p}} |u(\eta)|^p.\label{B}
\end{equation} 
\end{enumerate}
\end{cor}
\par\noindent{\bf Proof of Corollary \ref{cor2.3}}:
In Theorem \ref{T3.1} we  set $W_p(t)=w(t)^{p-1}= t^{\a p}$.  Then $ w(t)\in P(\mathbf R_+)$ for $\a \ge1/p'$,  
and $ w(t)\in \textcolor{black}{ Q}(\mathbf R_+)$ for   
$\a <1/p'$.
 For  $R>e$ we  set
\begin{equation}  \mu =
\begin{cases}& (\Lambda_{p}/ \Lambda_{\a,p})^{1/p} \eta^{(p-1-\a p)/(p-1)}, \quad \text{ if  }\a>1/p',
  \\ & \log R,  \qquad\quad \qquad \quad\quad \quad \qquad\quad \text{  if  }\a=1/p'.
\end{cases}
\end{equation}
When $\a \neq 1/p'$ we  have 
\begin{equation}
 F_\eta(t)= (\Lambda_{p}/ \Lambda_{\a,p})^{1/p} \, t  \quad\mbox{and}\quad  
 f_\eta(\eta)= \mu. 
\end{equation}
When $\a =1/p'$ we  have 
\begin{equation}
 F_\eta(t)=t \log (R\eta/t)\quad\mbox{and}\quad    f_\eta(\mu)= \mu.
\end{equation}
Then the  assertions (\ref{F}), (\ref{E}) and (\ref{B}) follow from  (\ref{nct1}).\qed

\begin{rem}\begin{enumerate}\item
The sharpness of  coefficients in  theses corollaries \ref{cor2.2}, \ref{cor2.3} will be  seen in Section 4 (4.2, Part 1).
\item If  $u\in C^1_c((0,\infty))$, the inequalities (\ref{G}), (\ref{F}) and  (\ref{B})   remain valid  for $\eta=+\infty$.
\end{enumerate}
\end{rem}
\par\medskip
In order to establish Hardy's inequalities in  a bounded domain of $\mathbf R^N$ as  an application, we need a  further refinement of 
 the previous results.
 
\begin{thm}\mbox{$(   w\in W_A(\mathbf R_+)= P_A(\mathbf R_+) \cup Q_A(\mathbf R_+))$ } \label{T3.2}\par \noindent
Assume that $1<p<\infty$, $\mu>0$ and  $ w\in W_A(\mathbf R_+)$. If $\eta>0$ is sufficiently small, then
 there exist  positive numbers    $C_0=C_0(w,p,\eta,\mu)$,  $C_1=C_1(w,p,\eta,\mu)$ and  $L=L(w,p,\eta,\mu)$  such  that  for  every $u\in W_{0}^{1,p}((0,\eta]; W_p)\cap C((0,\eta])$, we  have 
\begin{equation}
\begin{split}  \int_0^\eta  & \left( |u'(t)|^p -\frac{|u(t)|^p}{{F_\eta(t)}^p} \left(\Lambda_{p}    + \frac{C_0 }{G_\eta(t)^2}\right)\right) W_p(t)\,dt 
 \label{c1} \\ &\ge  C_1
 \int_0^\eta \left( |u'(t)|^p + \frac{|u(t)|^p}{{F_\eta(t)}^p} \left(\Lambda_{p} + \frac{C_0 }{G_\eta(t)^2}\right) \right) W_p(t)\, t\,dt +s(w)L |u(\eta)|^p,
\end{split}
\end{equation}
where $C_0$, $C_1$ and $L$ can be taken  independent of each $u$.

\end{thm}

Lastly we 
state a fundamental result which will  be  useful  in  the  subsequent.
\begin{prop}\label{ct1} $( w\in P(\mathbf R_+) )$\par\noindent
Assume that $ w\in P(\mathbf R_+)$ and $1<p<\infty$. Then  we have 
\begin{equation}\label{}
\inf_{u\in V_\eta}
\int_0^\eta  |u'(t)|^p W_p(t)\,dt =0,
 \end{equation}
 where $V_\eta= \{  u\in C^1([0,\eta]) :   u(0)=0, u(\eta)=1\}.  $
 \end{prop}
 \begin{rem}
The proof  will be  given in Section 6.
\end{rem}
\medskip
\subsection{Results in a  domain  of  $\mathbf R^N \,(N\ge 2)$ }

As an important application of one dimensional Hardy's inequalities in the  previous section, we describe Hardy's inequalities in a bounded domain 
$\Omega$ of  $\mathbf R^N \,(N\ge 2)$. Let $\delta(x):= { \rm dist}(x,\partial\Omega)$. 
For each  small $\eta >0$, $\Omega _\eta $ denotes a tubular neighborhood of $\partial \Omega$ and
   $\Sigma_\eta$ denotes the boundary of $ \Omega\setminus \Omega_\eta$, namely
 \begin{equation}
\Omega_\eta = \{ x\in \Omega:  \delta(x)<\eta \}\quad\mbox{and} \quad  \Sigma_\eta = \{ x\in \Omega:  \delta(x)=\eta \} .\label{NBD}
\end{equation} 
For the sake  of  simplicity,
by $W_p(\delta)$, $F_\eta(\delta)$ and $G_\eta(\delta)$ we denote 
$W_p(\delta(x))$, $F_\eta(\delta(x))$ and $G_\eta(\delta(x))$ respectively.\par
The  proofs of  Theorem \ref{CT3} and  Corollary \ref{NCC1}  will  be given in Section 5. Theorem \ref{NCT2} will be  proved in Section 6.

\begin{thm}\label{CT3} \mbox{\rm ( $   w\in W_A(\mathbf R_+)$ )} Assume  that  $\Omega$  is  a  bounded domain of  class $C^2$ in 
$\mathbf R^N$.
Assume that $1<p<\infty$ and   $ w\in W_A(\mathbf R_+)$. 
Assume  that $\mu>0$  and  $\eta $  is  a sufficiently  small positive  number. 
Then, there exist  positive  numbers    $C=C(w,p,\eta,\mu)$   and  $L'=L'(w,p,\eta,\mu)$ such  that  
for every  $ u\in W_0^{1,p}(\Omega;W_p(\delta))\cap C(\Omega)$, we have 
\begin{equation}\label{2.11}
\begin{split}
\int_{\Omega_\eta}   \left( |\nabla u|^p  -\Lambda_{p} \frac{|u|^p}{F_\eta(\delta)^p}  \right)W_p(\delta)\,dx
&\ge C \int_{\Omega_\eta}  \frac{|u|^p {W_p(\delta)}}{F_\eta(\delta)^p G_\eta(\delta)^2}\,dx \\
&+ s(w)L' \int_{\Sigma_\eta} |u|^p W_p(\delta)\,d\sigma_\eta,
 \end{split}
 \end{equation}
where $d\sigma_\eta$  denotes  surface elements on  $\Sigma_\eta$, and  $C$, $L'$ are independent of each $u$. 

\end{thm}
\begin{rem} We remark  that the  assumpsion $w\in W_A(\R_+)$ is needed 
even if  we do  not  have  the first term involving $G_\eta(t)$ in  the right-hand side. 
(See also Corollary \ref{NCC1}).
\end{rem}
\begin{cor}\label{NCC1}
Assume  that  $\Omega$  is  a  bounded domain of  class $C^2$ in $\mathbf R^N$.
Assume that $1<p<\infty$ and   $ w\in W_A(\mathbf R_+)$. Assume  that $\mu>0$ and  $\eta $  is  a sufficiently  small positive  number. 
Then, there exist  positive numbers  $\gamma=\gamma(w,p,\eta,\mu)$  and $ L'=L'(w,p,\eta,\mu)$ such  that  for every  
$ u\in W^{1,p}_{0}(\Omega;W_p(\delta))\cap C(\Omega)$, we have
\begin{equation}\label{2.7}
\int_{\Omega} \left( |\nabla u|^p  -\gamma\frac{|u|^p}{F_\eta(\delta)^p} \right)W_p(\delta)\,dx
\ge  s(w) L' \int_{\Sigma_\eta} |u|^pW_p(\delta)\,d\sigma_{\eta},
 \end{equation}
 where $d\sigma_\eta$  denotes  surface elements on  $\Sigma_\eta$, and  $\gamma$, $L'$ are independent of each $u$. 
 \end{cor}
Moreover 
 we have the followings:
 \begin{thm}\label{NCT2}\mbox{ \rm( $ w\in W_A(\mathbf R_+) $ )}  Assume  that  $\Omega$  is  a  bounded domain of   class $C^2$ in $\mathbf R^N$.
Assume that $1<p<\infty$, $\mu>0$ and   $ w\in W_A(\mathbf R_+)$. Then,  the followings are equivalent to each other.

\begin{enumerate}
\item  There exist  positive  numbers  $\gamma$, $\eta$  and $L'$  such  that the  inequality (\ref{2.7}) is valid for  every  $ u\in W_0^{1,p}(\Omega;W_p(\delta))\cap C(\Omega)$.

\item  For  a sufficiently  small  $\eta>0$,
there exist  positive  numbers  $\kappa$, $C$ and $L'$ such  that the  inequality (\ref{2.11})  with  $\Lambda_{p}$ replaced  by  $\kappa$ is  valid for   every  $ u\in W_0^{1,p}(\Omega;W_p(\delta))\cap C(\Omega)$.
\end{enumerate} 
\end{thm}
\par\medskip

\begin{prop}\label{CT1} \mbox{\rm ( $   w\in P(\mathbf R_+) $)} Assume  that  $\Omega$  is  a  bounded domain of  class $C^2$ in $\mathbf R^N$.
Assume that $1<p<\infty$ and   $ w\in P(\mathbf R_+)$. 
Then, for  an arbitrary  $\eta \in (0, \sup_{x\in \Omega} \delta(x) )$ we have 
\begin{equation}\label{2.8}
\inf
\left\{
\int_{\Omega}  |\nabla u|^p W_p(\delta) \,dx : u\in C^1_c (\Omega),  u=1 \mbox{  on } \{ \delta(x)=\eta\}\right \}=0.
 \end{equation}
 \end{prop}
 \begin{rem}\begin{enumerate}\item
 The proof  will be  given in Section 6 together with Proposition \ref{ct1}
 \item
From this we  see that constant functions belong to $W^{1,p}_{0}(\Omega;W_p(\delta))$, and  hence
 Hardy's inequality (\ref{HI}) never  holds when $w\in P(\mathbf R_+)$. 
\end{enumerate}
\end{rem}
\section{Proofs of Theorem \ref{T3.1}  and Theorem  \ref{T3.2} }
 
\subsection{  Lemmas  }

First we prepare  the following fundamental inequalities which are  established in \cite{anm}   as Lemma 2.1 for  $X>-1$. 
\begin{lem}\label{lem3.1} 
\begin{enumerate}\item
For $p\ge 2$ we have 
\begin{equation}
|1+X|^p-1-pX\ge  c(p)|X|^q, \quad \mbox{ for any } q\in [2,p]  \mbox{ and } X\in \mathbf {R}.
\label{3.0}\end{equation}
\item
For $1<p\le2$ and $M\ge1$, we  have 
\begin{equation} |1+X|^p-1-pX\ge  c(p)\begin{cases} &M^{p-2}X^2,\qquad \,\, |X|\le M,\\
& |X|^p,\qquad\qquad |X|\ge M.\end{cases}\label{3.1}\end{equation}
\end{enumerate}
Here $c(p)$ is a positive number independent of each $X$, $M\ge1$  and $q\in [2,p]$.
\end{lem}
{\noindent\bf Proof.} By  Taylar expansion  we have (\ref{3.0}) with $q=2$. 
For $p>1$, we  note that 
\begin{equation} \lim_{X\to 0} \frac{|1+X|^p-1-pX}{X^2}=\frac{p(p-1)}{2}, \quad   \lim_{|X|\to \infty} \frac{|1+X|^p-1-pX}{|X|^p} =1. \label{Add}
\end{equation}
Therefore  (\ref{3.0}) is valid for any $q\in [2,p]$ for a small $c(p)>0$.
If $X>-1$,  then (\ref{3.1}) also follows from Taylar expansion and  (\ref{Add}). 
If we choose $c(p)$ sufficiently small, then it
remains valid for $X\le -1$. 
\hfill$\Box$

\begin{rem}\label{density}
 $C_c^1((0,\eta]) $ is  densely  contained in $W_{0}^{1,p}((0,\eta]; W_p)$.
If $u\in  W_{0}^{1,p}((0,\eta]; W_p)$, then   $|u|\in W_{0}^{1,p}((0,\eta]; W_p)$ and 
$\|  u\|_{ W_{0}^{1,p}((0,\eta]; W_p)}$ $ =$ $ \| |u| \|_{ W_{0}^{1,p}((0,\eta]; W_p)}$.
Therefore, 
in  the  proofs of   Theorem \ref{T3.1}  and Theorem  \ref{T3.2}, without the loss of  generality  
 we  may assume that $u\in C_c^1((0,\eta]) $ and 
 $ u\ge 0$ in $(0,\eta)$.
 \end{rem}
For $u\in C_c^1((0,\eta])$ let us  set
\begin{equation}
u(t)= g_\eta(t) v(t),
\end{equation}
where $g_\eta(t)$ is  given by Definition \ref{fg}.
Clearly $v(t)\in C_c^1((0,\eta])$  and $u(t)'= g'_\eta(t) v(t)+g_\eta(t) v'(t)$. We define

\begin{equation} X(t)=  \begin{cases}&\frac{g_\eta(t)}{g'_\eta(t)} \frac{ v'(t)}{v(t)}=p'F_\eta(t) \frac{ v'(t)}{v(t)}
\textcolor{black}{s(w)} \quad \text{if } v(t)\neq 0,\\
& 0\,\,\quad\qquad\qquad\qquad\qquad\qquad\text{ if } v(t)= 0. \end{cases} \end{equation}
Then we have
\begin{equation}
u'(t)= g'_\eta(t) v(t)\left( 1+ X(t)\right), \quad \text{ if } v(t)\neq 0. \end{equation} 
Then it  follows   immediately from Lemma \ref{lem3.1} that  we have the following.
\begin{lem}\label{prop4.1}
Assume that  $\eta>0$, $\mu>0$, $X=X(t)$ and  $ w\in W(\mathbf R_+)$. 
Then, there exists a   positive number $c(p)$ such  that we  have  the followings:
\begin{enumerate}\item If $p\ge 2$,  then for any  $u\in C_c^1((0,\eta])$  and  $q\in [2,p]$ we  have
\begin{align}  \int_0^\eta   |u'(t)|^p W_p(t)\,dt
\ge \int_0^\eta |v(t)|^p |g_\eta'|^p  \left( 1+   pX +c(p)   |X|^q \right)W_p(t)\,dt. \label{1.lemma4.2}
\end{align}
\item If $1<p<2$, then for any  $u\in C_c^1((0,\eta])$ and $M\ge 1$ we  have
\begin{align}  \int_0^\eta   |u'(t)|^p W_p(t)\,dt
&\ge \int_0^\eta |v(t)|^p |g_\eta'|^p  \left( 1+   pX \right)W_p\,dt \label{2.lemma4.2} \\ 
+c(p) \int_0^\eta |v(t)|^p & |g_\eta'|^p    \left( M^{p-2}  X^2  \chi_{ \{|X|\le M\}} + |X|^p \chi_{\{|X|>M\}}\right)W_p(t)\,dt,\notag
\end{align}
where  $\chi_S(t)$ denotes  a characteristic function of a set $S$.
\end{enumerate}

\end{lem}
Now  we prepare the following lemma,
recalling a switching function;
\begin{equation} s(w)= \begin{cases} -1 \quad & \text{if } \quad w\in P(\mathbf R_+),\\
\,\,1\quad & \text{if } \quad w\in Q(\mathbf R_+).\end{cases}
\end{equation}
\begin{lem}\label{lemma4.2}Assume that  $\eta>0$, $\mu>0$ and  $ w\in W(\mathbf R_+)$.  Assume  that $u\in C_c^1((0,\eta])$, $X=X(t)$ and $ u(t)=g_\eta(t) v(t)$. Then  we  have:
\begin{enumerate}
\item
 \begin{equation}   |g'_\eta(t)|^{p-1} g_\eta(t) W_p(t)=1\label{4.13}
\end{equation}
\item
 \begin{equation}  p |v(t)|^p  |g_\eta'(t)|^p  X(t) W_p(t) =s(w) (|v(t)|^p)' \label{4.14}
\end{equation}

\item
 \begin{equation}  |v(t)|^p  |g_\eta'(t)|^p  W_p(t)= \Lambda_p \frac{ |u(t)|^pW_p(t)}{F_\eta(t)^p}\label{4.15}
\end{equation}

\item
 \begin{equation}|v(t)|^p  |g_\eta'(t)|^p  |X(t)|^2 W_p(t)= \frac{4p'}{p^2} | (v(t)^{p/2})'|^2 F_\eta(t) \label{4.16}
\end{equation}

\item
 \begin{equation}  |v(t)|^p  |g_\eta'(t)|^p  |X(t)|^p W_p(t)= (p')^{p-1} |v'(t)|^p F_\eta(t)^{p-1} \label{4.17}
\end{equation}
\end{enumerate}

\end{lem}
\noindent{\bf Proof of Lemma \ref{lemma4.2}:}
 \par\noindent
Proof of (\ref{4.13}):  By  Definition \ref{fg} (2) we see
 \begin{equation}  g'_\eta= (p')^{-1/p }  (f_\eta)^{-1/p } w^{-1} s(w). \label{4.18} \end{equation}  Then we have  (\ref{4.13}).
\par \noindent Proof of  (\ref{4.14}):  Using  (\ref{4.13}) together with
 \begin{equation}   p |v|^p  |g_\eta'|^p  X = s(w) (|v|^p)' |g'_\eta|^{p-1} g_\eta, \label{4.19}   \end{equation} 
 we have  (\ref{4.14}).
\par \noindent Proof of  (\ref{4.15}):  Noting  that  $p+p'=pp'$  and $(p')^{-p}= \Lambda_p$, we have 
 \begin{equation*} |v|^p|g_\eta'|^p=\frac{|u|^p}{ g^{p+p'}_\eta w^p}=\frac{|u|^p}{ g^{pp'}_\eta w^p}=\Lambda_p \frac{|u|^p}{ F_\eta^p}.
\label{4.20}   \end{equation*} 
\par \noindent Proof of  (\ref{4.16}):  Using (\ref{4.18}) and $2-(p-2)/(p-1)= p'$, we  have
 \begin{equation*}|v|^p  |g_\eta'|^p  |X|^2 W_p=|v|^{p-2} (v')^2 |g'_\eta|^{p-2}g_\eta^2 W_p= 
  \frac{4p'}{p^2} | (v^{p/2})'|^2 F_\eta. 
\end{equation*}
\par \noindent Proof of  (\ref{4.17}):  Using (\ref{4.18})  and $ g_\eta^p= (p' f_\eta)^{p-1}$ we  have 
\begin{equation*}  |v|^p  |g_\eta'|^p  |X|^p W_p= |v'|^p g_\eta^pW_p= (p')^{p-1} |v'|^p F_\eta^{p-1}.
\end{equation*}
\hfill $\Box$\par
By  the  definition of  $\Lambda_p$ we  have the following.
\begin{lem}\label{lem4.4}Assume that $\eta>0$  and $ w\in W(\mathbf R_+)$. Then we  have
\begin{equation}
\int_0^\eta (|v(t)|^p)' \,dt=\frac{ (\Lambda_p)^{1/p'}} {f_{\eta}(\eta)^{p-1}} |u(\eta)|^p, 
\label{4.21}
\end{equation}
where $v= u/g_\eta$ and $f_\eta(\eta)=\mu$, if $w\in P(\mathbf R_+)$; $\int_0^\eta w(s)^{-1}\,ds$, if $w\in Q(\mathbf R_+)$.

\end{lem}
By virtue of Lemma \ref{prop4.1}, Lemma \ref{lemma4.2} and Lemma \ref{lem4.4}  , we  have 
the  following:
\begin{lem}\label{prop4.2}
Assume that  $\eta>0$, $\mu>0$, $X=X(t)$ and  $ w\in W(\mathbf R_+)$. 
Then, there exists   a  positive number $c(p)$   such  that we  have  the followings:
 
\begin{enumerate}\item If $p\ge 2$,  then for any  $u\in C_c^1((0,\eta])$ we  have
\begin{align}  \int_0^\eta   |u'(t)|^p W_p(t)\,dt &\ge  \Lambda_p \int_0^\eta \frac{ |{u(t)}|^p W_p(t)\,dt }{F_\eta(t)^p} +s(w)\frac{ (\Lambda_p)^{1/p'}} {f_{\eta}(\eta)^{p-1}} |u(\eta)|^p \label{4.222}\\&
+ d(p)\int_0^\eta |(|v(t)|^{p/2})'|^2 F_\eta(t) \,dt,\notag
\end{align}
where  $d(p)= c(p) 4p'/p^2$.\par
Moreover the last term can be  replaced by  $c(p)(p')^{p-1}\int_0^\eta |v'(t)|^p F_\eta(t)^{p-1} \,dt$.
\item If $1<p<2$, then for any  $u\in C_c^1((0,\eta])$ and $M\ge 1$ we  have
\begin{align} & \int_0^\eta   |u'|^p W_p(t)\,dt
\ge \Lambda_p \int_0^\eta \frac{ |{u(t)}|^p W_p(t)\,dt }{F_\eta(t)^p} + s(w)\frac{ (\Lambda_p)^{1/p'}} {f_{\eta}(\eta)^{p-1}} |u(\eta)|^p  \label{4.23}\\ &\notag
+c(p)(p')^{p-1}\int_0^\eta |v'(t)|^p F_\eta(t)^{p-1}  \chi_{\{|X|>M\}}\,dt \\&+ M^{p-2}d(p)\int_0^\eta
|(|v(t)|^{p/2})'|^2 F_\eta(t)  \chi_{ \{|X|\le M\}} \,dt, \notag
\end{align}
where  $d(p)= c(p) 4p'/p^2$ and $\chi_S(t)$ denotes  a characteristic function of a set $S$.
\end{enumerate}
\end{lem}
\par\noindent{\bf Proof of Lemma \ref{prop4.2}:}
\par\noindent{Proof of (\ref{4.222}): } This follows from (\ref{1.lemma4.2}) of
 Lemma \ref{prop4.1}, (\ref{4.13}), (\ref{4.14}), (\ref{4.15}), (\ref{4.16}) of Lemma \ref{lemma4.2} and Lemma \ref{lem4.4}. By using (\ref{4.17}) instead of (\ref{4.16}), the alternative inequality follows. 
 \par\noindent{Proof of (\ref{4.23}): } This follows from (\ref{2.lemma4.2}) of
 Lemma \ref{prop4.1}, (\ref{4.13}), (\ref{4.14}), (\ref{4.15}), (\ref{4.16}), (\ref{4.17}) of Lemma \ref{lemma4.2} and Lemma \ref{lem4.4}. \qed
\begin{rem}\label{rem4.2}
 From
 (\ref{4.222}) with the last term 
 replaced by  $c(p)(p')^{p-1}\int_0^\eta |v'|^p F_\eta^{p-1} \,dt$, we   have the following  particular inequality for $p\ge 2$:

\begin{align}  \int_0^\eta  & |u'|^p W_p\,dt 
\ge  \Lambda_p \int_0^\eta \frac{ |{u(t)}|^p W_p(t)\,dt }{F_\eta(t)^p} +s(w)\frac{ (\Lambda_p)^{1/p'}} {f_{\eta}(\eta)^{p-1}} |u(\eta)|^p\\ &
+\frac{d(p)}{2}\int_0^\eta |(|v(t)|^{p/2})'|^2 F_\eta(t) \,dt +\frac{c(p)(p')^{p-1}}{2}\int_0^\eta |v'(t)|^p F_\eta(t)^{p-1} \,dt.\notag
\end{align}

\end{rem}

\subsection{Proof  of  Theorem \ref{T3.1}}
\par\noindent{\bf Part 1: Proof of the inequality (\ref{nct1})  and its sharpness }\par
Assume  that  $w\in W(\mathbf R_+)= P(\mathbf R_+)\cup Q(\mathbf R_+)$.
It follows  from Lemma \ref{prop4.2},
we  clearly have the inequality (\ref{nct1})  in Theorem \ref{T3.1}. Hence  we proceed to  the proof of optimality of the coefficients $\Lambda_p$  and $(\Lambda_p)^{1/p'}$ in  the right-hand side of (\ref{nct1}). By  the density argument,
 we can adopt as test functions
\begin{align} u_\varepsilon(t) = f_\eta(t)^{1/p'+s(w)\varepsilon } \qquad (0<\varepsilon).
\end{align}
We  note  that
\begin{equation}  u'_\varepsilon(t)= \left( 1/p' +s(w)\varepsilon\right) f_\eta(t)^{s(w)\varepsilon -1/p}
{s(w)}{w}^{-1}.
\end{equation}
Then the left-hand side becomes
\begin{align*} 
\int_0^\eta|u_\varepsilon'(t)|^pW_p(t)\,dt&= \left( 1/p' +s(w)\varepsilon\right)^p\int_0^\eta 
 f_\eta(t)^{s(w)\varepsilon p-1}{w(t)}^{-1}\,dt \notag\\
&= \left( 1/p' +s(w)\varepsilon\right)^p {f_\eta(\eta)^{s(w) \varepsilon p}}({p\varepsilon})^{-1}.
\end{align*}
In  a similar way the right-hand side becomes
\begin{align*}\Lambda_p &\int_0^\eta  \frac{|u_\varepsilon(t)|^p W_p(t)}{F_\eta(t)^p }\,dt+ s(w) \Lambda_p^{1/p'} \frac{|u_\varepsilon (\eta)|^p}{f_\eta(\eta) ^{p-1}}\notag\\
&= \left( 1/{p' }\right)^{p}\int_0^\eta  f_\eta(t)^{s(w)\varepsilon p-1}w(t)^{-1}\,dt+ s(w) \Lambda_p^{1/p'} \frac{|u_\varepsilon (\eta)|^p}{f_\eta(\eta) ^{p-1}}
 \notag\\
&= \left( 1/p' \right)^{p} {f_\eta(\eta)^{s(w) \varepsilon p}}({p\varepsilon})^{-1} + s(w) \Lambda_p^{1/p'} \frac{|u_\varepsilon (\eta)|^p}{f_\eta(\eta) ^{p-1}}.
\end{align*}
Finally we  reach to
\begin{equation} \left(\left( 1/p' +s(w)\varepsilon\right)^p-\left( 1/{p' } \right)^{p}\right)
\ge s(w) \Lambda_p^{1/p'}p\varepsilon \qquad (0<\varepsilon ),
\end{equation}
and this inequality is  clear from the    convexity. \hfill $\Box$
\par\medskip
\par\noindent{\bf Part 2: Proof of the inequality (\ref{nct2}) }\par
In order to establish  the  inequality (\ref{nct2})  in Theorem \ref{T3.1} we  shall estimate  the term involving 
$ |(|v|^{p/2})'|^2 F_\eta^{}$ from  below  and   
 use a positive term involving $ |v'|^p F_\eta^{p-1}$ to absorb negative error terms.
To  this end we need more  notations:\par
For $u,v \in C^1_c((0,\eta])$, we  retain 
$$
u(t)= g_\eta(t) v (t)\quad \mbox{and }\quad X(t)=  \frac{g_\eta(t)}{g'_\eta(t)} \frac{ v'(t)}{v(t)}=p'F_\eta(t) \frac{ v'(t)}{v(t)}
s(w)\quad \textcolor{black}{(v(t)\neq0)}.
$$
By Remark \ref{density}, it  suffices to assume in the proof  that $u$ and $v$ belong  to the class 
$ \textcolor{black}{G((0,\eta])}$ defined by the following:
\begin{df}\label{df4.2} 
\begin{equation}\label{4.22} G((0,\eta])=\{ v\in C_c^1((0,\eta]):    v(t)\ge 0 \,\text { in } (0,\eta]\, \} \end{equation}
\end{df}  
\begin{df}\label{df4.3} For $v \in G((0,\eta])$ and $M>1$ we define  two subsets of $[0,\eta]$ as
follows:
 \begin{equation} \begin{cases}& A(v,M)= \{ t\in [0,\eta] : \,\, \,
 p' F_\eta(t) |v'(t)| \le M {v(t)}\}\\
& B(v,M)= \{ t\in [0,\eta] :\,\, \,  p' F_\eta(t) |v'(t)| > M {v(t)}\}\\
\end{cases}
\end{equation}
$A(v,M)$  and  $B(v,M)$ are sometimes abbreviated as $A$ and $B$ respectively.
\end{df}
\begin{rem}
 From  Definition \ref{df4.2}  and Definition \ref{df4.3} we  see $A\cup B=[0,\eta]$. 
\textcolor{black}{ If $v(t)\in G((0,\eta])$ and $v(s)=0$ for some $s\in [0,\eta)$, then  $v'(s)=0$, and hence $s\in A(v,M)$ for any $M>0$. }
 We note  that the  set 
$C(v,M)=\{ t\in [0,\eta] :\,\, \,  p' F_\eta(t) |v'(t)| = M {|v(t)|}\}$ coincides with the  set of critical points of $\, \Psi(t):=\log (|v|f_\eta^{\pm M/p'})$. Namely, $\{ t\in[0,\eta] : \Psi'(t)=0 \}$.
\end{rem}
The  following lemma will be  established in Appendix. 
\begin{lem}\label{lemma4.4}Assume that  $\eta>0$, $\mu>0$, $M>1$ and  $ w\in W(\mathbf R_+)$. For  any  
$v\in  G((0,\eta])$ we have the followings: 
\begin{enumerate}\item For $A=A(v,M)$ and $B=B(v,M)$,
\begin{equation}\label{4.32}
\begin{split} &\int_A  |(v(t)^{p/2})'|^2  F_\eta(t)\,dt \ge -\frac{ (\Lambda_p)^{1/p'}} {2\mu f_{\eta}(\eta)^{p-1}} |u(\eta)|^p
+\frac{1}{4(p')^{p-1}}\int_A\frac{|u(t)|^pW_p(t)}{F_\eta(t)^pG_\eta(t)^2} \,dt  \\
&+\frac{1}{2(p')^{p-1}}\int_B\frac{|u(t)|^pW_p(t)}{F_\eta(t)^pG_\eta(t)^2}\,dt- \frac{p}{2\mu}\left(\frac{p'}{M}\right)^{p-1}\int_B|v'(t)|^p F_\eta(t)^{p-1}\,dt.
\end{split}
\end{equation}
\item 
\begin{equation}\label{4.32'} 
\begin{split}
\int_0^\eta  |(v(t)^{p/2})'|^2  F_\eta(t)\,dt &\ge -\frac{ (\Lambda_p)^{1/p'}} {2\mu f_{\eta}(\eta)^{p-1}} |u(\eta)|^p\\&
+\frac{1}{4(p')^{p-1}}\int_0^\eta \frac{|u(t)|^pW_p(t)}{F_\eta(t)^pG_\eta(t)^2} \,dt 
\end{split}
\end{equation}

\end{enumerate}
\end{lem}
\begin{rem}
The assertion 2 follows from Lemma \ref{lemma8.1} with $S=[0,\eta]$.
 The inequality (\ref{nct2})  will follow from  this  lemma together with Lemma \ref{prop4.2}.
 \end{rem}
\par\noindent
{\bf End  of  the proof of (\ref{nct2}):}
If $p\ge 2$,  then for any  $u\in G((0,\eta])$ we  have by using  (\ref{4.32'}) of  Lemma \ref{lemma4.4} 
\begin{align*}  \int_0^\eta   |u'(t)|^p W_p(t)\,dt
&\ge s(w)\frac{ (\Lambda_p)^{1/p'}} {f_{\eta}(\eta)^{p-1}} |u(\eta)|^p+ \Lambda_p \int_0^\eta \frac{ |{u(t)}|^p W_p(t)\,dt }{F_\eta(t)^p}\\ &
+d(p)\left(-\frac{ (\Lambda_p)^{1/p'}} {2\mu f_{\eta}(\eta)^{p-1}} |u(\eta)|^p
+\frac{1}{4(p')^{p-1}}\int_0^\eta\frac{|u(t)|^pW_p(t)}{F_\eta(t)^pG_\eta(t)^2}\,dt\right)
\notag\\
&= s(w)\frac{ (\Lambda_p)^{1/p'}} {f_{\eta}(\eta)^{p-1}} \left( 1- s(w)\frac{d(p)}{2\mu}\right) |u(\eta)|^p +\Lambda_p \int_0^\eta \frac{ |{u(t)}|^p W_p(t)\,dt }{F_\eta(t)^p}\,dt\notag
\\ &+\frac{d(p)}{4(p')^{p-1}}\int_0^\eta \frac{|u(t)|^pW_p(t)}{F_\eta(t)^pG_\eta(t)^2}\,dt.\notag
\end{align*}
 Here  we assume that $d(p) $ is  so small that $ d(p)/\mu<1$ if $w\in Q(\mathbf R_+)$.
Then we  get (\ref{nct2}) with 
 \begin{equation}\begin{cases}&L= { (\Lambda_p)^{1/p'}} {f_{\eta}(\eta)^{1-p}} \left( 1-s(w) {d(p)}/({2\mu})\right),\\  &\\&
 C={d(p) (\Lambda_p)^{1/p'}}/{4}.\end{cases}
 \end{equation}

 If $1<p<2$, then for any  $u\in G((0,\eta])$ and $M\ge 1$ we  have
\begin{align*} & \int_0^\eta   |u'|^p W_p\,dt
\ge s(w)\frac{ (\Lambda_p)^{1/p'}} {f_{\eta}(\eta)^{p-1}} |u(\eta)|^p+ \Lambda_p \int_0^\eta \frac{ |{u(t)}|^p W_p(t)\,dt }{F_\eta(t)^p}\\ &
+c(p)(p')^{p-1}\int_{B(v,M)} |v'(t)|^p F_\eta(t)^{p-1} \,dt+ M^{p-2}d(p)\int_{A(v,M)}
|(|v(t)|^{p/2})'|^2 F_\eta  \,dt. \notag
\end{align*}
Here  $c(p)$ is  a positive number independent of each $u$ and $d(p)= c(p) 4p'/p^2$.
For  the last term we use (\ref{4.32}) of Lemma \ref{lemma4.4},  then we have
\begin{align*} & \int_0^\eta   |u'(t)|^p W_p(t)\,dt
\ge 
\frac{ s(w)(\Lambda_p)^{1/p'}} {f_{\eta}(\eta)^{p-1}} \left( 1-s(w) \frac{d(p)M^{p-2}}{2\mu}\right) |u(\eta)|^p \\&+\Lambda_p \int_0^\eta \frac{ |{u(t)}|^p W_p(t)\,dt }{F_\eta(t)^p}
+c(p)(p')^{p-1}\left( 1- \frac{2}{\mu(p-1)M} \right)
\int_{B(v,M)} |v'(t)|^p F_\eta(t)^{p-1} \,dt 
\\&
+\frac{d(p)M^{p-2} (\Lambda_p)^{1/p'}}{2}\left( \frac12\int_{A(v,M)} \frac{|u(t)|^pW_p(t)}{F_\eta(t)^pG_\eta(t)^2}\,dt+
\int_{B(v,M)}  \frac{|u(t)|^pW_p(t)}{F_\eta(t)^pG_\eta(t)^2}\,dt\right).
\notag
\end{align*}

Then we  take a sufficiently large  $M$ so    that  we  have $1- {2}/({\mu(p-1)M}) >0$ and ${d(p)M^{p-2} }/({2\mu})<1$ if $w\in Q(\mathbf R_+)$. Then
\begin{align*} 
 &\int_0^\eta   |u'(t)|^p W_p(t)\,dt
\ge 
 s(w)L |u(\eta)|^p +\Lambda_p \int_0^\eta \frac{ |{u(t)}|^p W_p(t)\,dt }{F_\eta(t)^p}
+C \int_0^\eta \frac{|u(t)|^pW_p(t)}{F_\eta(t)^pG_\eta(t)^2}\,dt,
\end{align*}
where \begin{equation}\begin{cases}&L= { (\Lambda_p)^{1/p'}} {f_{\eta}(\eta)^{1-p}} \left( 1-s(w) {d(p)M^{p-2}}/({2\mu})\right),\\  &\\&
 C={d(p)M^{p-2} (\Lambda_p)^{1/p'}}/{4}.\end{cases}
 \end{equation}
\qed
\subsection{Proof  of  Theorem \ref{T3.2}}
Let $v\in G((0,\eta])$.
Recall that
$
u(t)= g_\eta(t) v(t)$, $u'(t)= g'_\eta(t) v(t)\left( 1+ X(t)\right) $ $(v\neq0)$ and 
\begin{equation*} X(t)=  \frac{g_\eta(t)}{g'_\eta(t)} \frac{ v'(t)}{v(t)}=p'F_\eta (t)\frac{ v'(t)}{v(t)}\textcolor{black}{s(w)}\quad (v(t)\neq 0); 0\quad (v(t)=0).
\end{equation*}
Then we have two elementary lemmas.
\begin{lem}\label{lemma4.7} For  any $X\in \mathbf R$  and  any $M>1$,
\begin{equation} |1+X|^p \le \begin{cases}& (1+M)^p,  \qquad  |X|\le M, \\
& 2^p |X|^p,\qquad \quad |X|> M.
 \end{cases}
\end{equation}
\end{lem}
\begin{lem}\label{lemma4.8}For any $M>1$,
\begin{equation}  |u'(t)|^p W_p(t)\le  \begin{cases}&  \Lambda_p(1+M)^p  { u(t)^p W_p(t)}{ F_\eta(t)^{-p}}, \quad\qquad |X(t)|\le M,\\
& 2^p\Lambda_p^{-1/p'} |v'(t)|^p F_\eta(t)^{p-1},\qquad \qquad \quad \, \,\, |X(t)|> M.
 \end{cases}
\end{equation}
\end{lem}
\par\noindent{\bf Proof:} Since $u'(t)= g'_\eta(t) v(t)\left( 1+ X(t)\right) $, we get 
$$|u'(t)|^pW_p(t)\le |g'_\eta(t)|^pv(t)^p|\,1+X(t)|^p \textcolor{black}{W_p(t).}$$
 If $|X(t)|\le M$,  then from Lemma \ref{lemma4.7} and (\ref{4.15})
we  get the desired estimate. If $|X(t)|\ge M$,  then from Lemma \ref{lemma4.7} and (\ref{4.17}) we get  the desired one.\hfill $\Box$\par
By \textcolor{black}{ Proposition \ref{prop2.1} } we have  the following.
\begin{lem}\label{lemma4.9} Assume  that $w$ is admissible in  the  sense of Definition \ref{admissibility}. Then, there are some  positive  numbers $K$  and $\eta$  such  that we have for  
any $u\in C_c^1(0,\eta])$
\begin{equation}\int_0^\eta  \frac{ u(t)^pW_p(t)}{ F_\eta(t)^p}t\,dt \le K^2\int_0^\eta  \frac{ |u(t)|^p W_p(t) \,dt}{ F_\eta(t)^p G_\eta(t)^2}.\end{equation}
\end{lem}
\par\noindent{\bf Proof of Theorem \ref{T3.2}:}
For a sufficiently small $\eta>0$, temporally we  set 
\textcolor{black}{\begin{equation} Q_{ L}(u)= \int_0^\eta   |u'(t)|^p W_p(t)\,dt 
  -\Lambda_p \int_0^\eta \frac{ |{u(t)}|^p W_p(t)\,dt }{F_\eta(t)^p}-s(w) {L}|u(\eta)|^p,
\end{equation}}
where $L$ is a positive constant independent of each $u$, \textcolor{black}{which will be specified later}.
In order to prove Theorem \ref{T3.2}, it  suffices to control three terms below by $Q_{L}(u)$  :
$$\int_0^\eta \frac{ |u(t)|^p W_p(t) \,dt}{ F_\eta(t)^p G_\eta(t)^2},  \quad \int_0^\eta  \frac{ u(t)^pW_p(t)}{ F_\eta(t)^p}\, t\,dt \quad\text{and} \quad \int_0^\eta |u'(t)|^p W_p(t) \, t\,dt.$$
It follows  form Theorem \ref{T3.1} (\ref{nct2}) that we  have for $D=C^{-1}$
$$\int_0^\eta \frac{ |u(t)|^p W_p(t) \,dt}{ F_\eta(t)^p G_\eta(t)^2}\le D\cdot Q_L(u),$$
where $L=L(w,p,\eta,\mu)$ is the same constant  in (\ref{nct2}).
From this and Lemma \ref{lemma4.9}   we  also have 
$$\int_0^\eta  \frac{ u(t)^pW_p(t)}{ F_\eta(t)^p}\, t\,dt \le K^2 D\cdot Q_L(u).$$
By using Lemma \ref{lemma4.8} we  see that
\begin{align*}\int_0^\eta |u'(t)|^p W_p(t) \, t\,dt &\le \Lambda_p(1+M)^p\int_{|X|\le M} \frac{ u(t)^pW_p(t)}{ F_\eta(t)^p}\, t\,dt \\&+2^p\Lambda_p^{1/p-1}\int_{|X|> M}|v'(t)|^p F_\eta(t)^{p-1}\, t\,dt.\end{align*}
Then,  the first term in the right-hand side can be controlled by  the  same $Q_L(u)$.
Now  we claim that   for  some  positive number $K'$ \begin{equation} \int_{|X|> M}|v'(t)|^p F_\eta(t)^{p-1}t\,dt\le K'\cdot Q_{L}(u), \label{4.43}\end{equation}
where  $L= {(\Lambda_p)^{1/p'}}{f_\eta(\eta)^{1-p}} $.
\par\noindent
First  we assume  that $p\ge 2$. 
By  (\ref{4.222})  of  Lemma \ref{prop4.2} with the last term replaced  by $$c(p)(p')^{p-1}\int_0^\eta |v'(t)|^p F_\eta(t)^{p-1} \,dt,$$we  see that  (\ref{4.43}) is valid, provided that  $\eta$ is  sufficiently small (See also Remark \ref{rem4.2}). Secondly  we  assume  that  $1<p<2$. 
Then in a similar way  the assertion follows from  (\ref{4.23}) of Lemma \ref{prop4.2}.
As a result we have the desired inequality (\ref{c1}) for some positive numbers $C_0, C_1$ and $L$ which are independent of each $u$.
\hfill$\Box$
\section{ Proofs  of  Theorem \ref{CT3}  and  Corollary \ref{NCC1}}
We  establish Theorem \ref{CT3}  and Corollary  \ref{NCC1} using Theorem \ref{T3.2}.
\par \medskip
\noindent{\bf  Proof of Theorem \ref{CT3}:} 
Let  us  prepare some notations and  fundamental facts.
Define $\Sigma = \partial \Omega$  and $\Sigma_t =\{ x\in  \Omega : \delta(x)= t\}$.
Since   $\Sigma$ is  is  of  class $C^2$, there exists an $\eta_0>0$  such  that we  have 
a $C^2$ diffeomorphism $ G: \Omega_\eta \mapsto (0,\eta)\times \Sigma$
for any  $\eta\in (0,\eta_0)$.  By $G^{-1}(t, \sigma) \, ((t,\sigma) \in (0,\eta_0) \times \Sigma)$ we denote the inverse of $G$. Let $H^{}_t$ denote the mapping 
$G^{-1}(t, \cdot) $ of $\Sigma$ onto $\Sigma_t$.
This mapping is also a $C^2$ diffeomorphism  and its Jacobian is close to $1$ in $(0,\eta_0)\times \Sigma$.
Therefore,
for every non-negative  continuous function $u$ on $\overline{\Omega_\eta}$ with $\eta\in (0,\eta_0)$  we  have 
\begin{align}& \int_{\Omega_\eta}u = \int _0^\eta \,dt \int_{\Sigma_t} u \,d\sigma_t =\int _0^\eta \,dt \int_\Sigma u( H_t(\sigma)) ( \rm{Jac}\, H_t )\,d\sigma, \label{5.1}\\
& | \rm{Jac}\, H_t (\sigma)-1|\le ct, \quad \mbox{  for every   } (t,\sigma)\in (0,\eta_0)\times  \Sigma, \label{5.2}
\end{align}
where $c$  is  a positive  constant independent of  each $(t,\sigma)$, $d\sigma$  and $d\sigma_t$  denote  surface elements on $\Sigma$  and $\Sigma_t$  respectively.
Then we  have 
 \begin{align} 
 & 
\int_\Sigma u( H_\eta(\sigma)) ( 1-c\eta)\,d\sigma \le 
\int_{\Sigma_\eta}u\,d\sigma_\eta\le   \int_\Sigma u( H_\eta(\sigma)) ( 1+c\eta)\,d\sigma. \label{4.2}
\end{align}
Again for the sake  of  simplicity, we denote  $W_p(\delta(x))$ and $F_\eta(\delta(x))$ by $W_p(\delta)$ and $F_\eta(\delta)$ respectively, and
the symbol $dx$ is often  abbreviated.
Then we immediately  have for $v=u(H_t(\sigma))$
\begin{align}
   \int_{\Sigma} \,d\sigma \int_0^\eta \left |\frac{\partial v}{\partial t} \right |^p (1-ct) &
  W_p(t) \,dt \le \int_{\Omega_\eta}|\nabla u|^p W_p(\delta)\,dx, \notag\\
 \int_{\Sigma} \,d\sigma \int_0^\eta \frac{| v|^pW_p(t) }{F_\eta(t)^p} (1-ct)  \,dt &\le  \int_{\Omega_\eta}\frac{| u|^pW_p(\delta)\,dx }{F_\eta(\delta)^p}\le  \int_{\Sigma} \,d\sigma \int_0^\eta\frac{| v|^pW_p(t) }{F_\eta(t)^p} (1+ct) \,dt. \notag
 \end{align}
 \noindent{\bf  Proof  of (\ref{2.11}):} Under these consideration, 
 (\ref{2.11}) is  reduced  to  the following one dimensional Hardy's inequality. Setting $v(t) =u(H_t(\sigma))$  and $v'= \partial v/\partial t$ we  have
 \begin{equation}
 \begin{split}  \int_0^\eta  &  |v'(t)|^p W_p(t)(1-ct)\,dt \ge
 \Lambda_{p}   \int_0^\eta \frac{|v(t)|^p W_p(t)}{F_\eta(t)^p} (1+ct)\,dt \\ 
 &+ C \int_0^\eta \frac{|v(t)|^p W_p(t)}{F_\eta(t)^p \textcolor{black}{G_\eta(t)^2} }(1+ct)\,dt  
+ s(w) L' |v(\eta)|^p W_p(\eta)(1+s(w)c\eta).
\end{split}
\end{equation}
Equivalently we have
 \begin{equation}  \label{5.7} 
 \begin{split} \int_0^\eta  & \left( |v'(t)|^p -\frac{|v(t)|^p}{F_\eta(t)^p}  \left(\Lambda_{p}    + \frac{C }{G_\eta(t)^2}\right)\right) W_p(t)\,dt 
\\ &\ge  c
 \int_0^\eta \left( |v'(t)|^p + \frac{|v(t)|^p}{F_\eta(t)^p} \left(\Lambda_{p} + \frac{C }{G_\eta(t)^2}\right) \right) W_p(t)\, t \,dt \\&
 +s(w) L' |v(\eta)|^p W_p(\eta)(1+s(w)c\eta).
 \end{split}
\end{equation}
Assume that $0<\eta\le  1/(2c)$. Then $1-c\eta\ge 1/2$ and   (\ref{5.7}) clearly  follows  from Theorem \ref{T3.2} with 
 $ C_0=C $, $C_1=c$    and  $L =L'W_p(\eta)(1+s(w)c\eta)$. \hfill $\Box$
 \par\medskip
\noindent{\bf  Proof of Corollary \ref{NCC1}: } 
First we  treat the  case  that \textcolor{black}{$w\in P_A(\mathbf R_+)$}. 
Assume that Hardy's inequality (\ref{2.7}) does  not  hold. 
Then there exists a sequence of  functions 
$\{u_k\} \subset W^{1,p}_{0}(\Omega; W_p(\delta)) \cap C(\Omega)$ 
such  that
\begin{equation} \begin{cases}&\lim_{k\to \infty}\left(\int_{\Omega} |\nabla u_k|^pW_p(\delta)\,dx +   \int_{\Sigma_\eta}|u_k|^p W_p(\delta)\,d\sigma_\eta\right)=0,\\ &
 \int_{\Omega} {|u_k|^pW_p(\delta)}/{F_\eta^p(\delta)}\, dx =1 \quad (k=1,2,\cdots). 
 \end{cases} \label{5.8}
\end{equation}
For  a sufficiently small $\eta>0$,
let  $W^{1,p}(\Omega\setminus \overline{\Omega_\eta}; W_p(\delta))$  be 
 given by the  completion of $C^\infty(\Omega\setminus \overline{\Omega_\eta} )$ with respect to the norm 
defined by
\begin{equation*} 
\| u\|_{W^{1,p}(\Omega\setminus \overline{\Omega_\eta}; W_p(\delta)) } =
\| |\nabla u|  \|_{ L^p(\Omega\setminus \overline{\Omega_\eta}; W_p(\delta))} + \| u\|_{ L^p(\Omega\setminus \overline{\Omega_\eta}; W_p(\delta))}. 
\end{equation*}
Since $W_p(\delta)>0$ in $\overline{\Omega\setminus{\Omega_\eta}}$, $W^{1,p}(\Omega\setminus \overline{\Omega_\eta}; W_p(\delta)) $  is well-defined and  becomes
a Banach space  with the norm $\| \cdot\|_{W^{1,p}(\Omega\setminus \overline{\Omega_\eta}; W_p(\delta)) }
$.
We  note  that
$\int_{\Sigma_\eta}|u_k|^p W_p(\delta)\,d\sigma_\eta$ $(k=1,2,3\ldots)$ is  bounded, because  the  trace operator $T: W^{1,p}(\Omega\setminus \overline{\Omega_\eta}; W_p(\delta))$  $\mapsto $ 
 $L^p(\Sigma_\eta; W_p(\eta) )$ is  continuous.
By Theorem {\ref{CT3}} we  have 
\begin{align*} 
&\int_{\Omega} |\nabla u_k|^pW_p(\delta)\,dx =  \int_{\Omega_\eta} |\nabla u_k|^pW_p(\delta) \,dx
+ \int_{\Omega\setminus \Omega_\eta} |\nabla u_k|^pW_p(\delta)\,dx \\
&\ge  \Lambda_{p}\left( 1-  \int_{\Omega \setminus \Omega_\eta}\frac{ {|u_k|^p W_p(\delta)}}{F_\eta^p(\delta)}\,dx\, \right)
+ \int_{\Omega\setminus \Omega_\eta} |\nabla u_k|^pW_p(\delta)\,dx-L' \int_{\Sigma_\eta} |u_k|^pW_p(\delta)
\,d\sigma_\eta. 
\end{align*}
Since  $\delta \ge \eta$ in $\Omega\setminus \Omega_\eta$,  by  the  standard argument we have $u_k\to C \,( constant)$ in $W^{1,p}(\Omega\setminus \overline{\Omega_\eta}; W_p(\delta))$ as $k\to \infty$. By (\ref{5.8}) we  have $C=0$.
Hence  we see
$ 0 \ge \Lambda_{p}$, and   we  reach  to   a contradiction. 
\par Secondly we  treat  the case that $w\in Q(\mathbf R_+)$. In  this  case  we assume that
\begin{align} &\lim_{k\to \infty}\int_{\Omega} |\nabla u_k|^pW_p(\delta)\,dx=0, \quad
 \int_{\Omega} \frac{|u_k|^pW_p(\delta)}{F_\eta^p(\delta)}\,dx =1 \quad (k=1,2,\cdots). \label{5.9}
\end{align}
Then, 
by Theorem {\ref{CT3}} we  have 
\begin{align} 
&\int_{\Omega} |\nabla u_k|^pW_p(\delta)\,dx =  \int_{\Omega_\eta} |\nabla u_k|^pW_p(\delta)\,dx + \int_{\Omega\setminus \Omega_\eta} |\nabla u_k|^pW_p(\delta) \,dx \label{5.10}\\
&\ge  \Lambda_{p}\left( 1-  \int_{\Omega \setminus \Omega_\eta}\frac{ {|u_k|^p W_p(\delta)}}{F_\eta^p(\delta)}\,dx \right)
+ \int_{\Omega\setminus \Omega_\eta} |\nabla u_k|^pW_p(\delta)\,dx+L' \int_{\Sigma_\eta} |u_k|^pW_p(\delta)\,d\sigma_\eta. \notag
\end{align}
Again we have $u_k\to C \,( constant)$ in $W^{1,p}(\Omega\setminus \overline{\Omega_\eta}; W_p(\delta))$ as $k\to \infty$. By (\ref{5.9}) and (\ref{5.10}) we  have $C=0$.
Hence  
 we see
$ 0 \ge \Lambda_{p}$, and  this  is   a contradiction.

 \hfill $\Box$

\section{Proof of Theorem \ref{NCT2} } 
By virtue of  Corollary \ref{NCC1}  and  its  proof, it suffices  to  show  the  implication $1\to 2$. Since $F_\eta^{-1}\notin L^1((0,\eta))$, for an  arbitrary $\varepsilon>0$ we have 
$G_\eta (\delta)^{-1} <\varepsilon$ in $ \Omega_\eta$ provided  that $\eta$  is sufficiently  small. 
Therefore we assume that $C=0$ without the loss of generality.
\par\medskip
Here we introduce a Banach space ${W_0^{1,p}(\Omega_\eta; W_p(\delta)) }$ which corresponds  to $W_{0}^{1,p}((0,\eta]; W_p)$  defined in Subsection 3.1. 
Let $ \tilde{C}_c^\infty(\Omega_\eta)= \{  u|_{\Omega_\eta}: u\in C_c^\infty(\Omega)\}$, where
by $ u|_{\Omega_\eta}$ we  denote the restriction of $u $ to $\Omega_\eta$. Let  $W_0^{1,p}(\Omega_\eta; W_p(\delta))$ be 
 given by the  completion of $\tilde{C}_c^\infty(\Omega_\eta)$ with respect to the norm 
defined by
\begin{equation*} 
\| u\|_{W_0^{1,p}(\Omega_\eta; W_p(\delta)) } =
\| |\nabla u|  \|_{ L^p(\Omega_\eta; W_p(\delta))} + \| u\|_{ L^p( \Omega_\eta; W_p(\delta))}. 
\end{equation*}
Then, 
 $W_0^{1,p}(\Omega_\eta; W_p(\delta)) $  is a Banach space  with the norm $\| \cdot\|_{W^{1,p}_{0}(\Omega_\eta; W_p(\delta)) }
$.
Then we  prepare a lemma on   extension:
\begin{lem}\label{l11}\mbox{\rm ( Extension )} Assume that $\eta_0$ is a sufficiently small positive  number. Then
for  any $\eta \in (0,\eta_0) $ there exists an extension operator  
$E=E(\eta) : W_0^{1,p}(\Omega_\eta; W_p(\delta))\mapsto  W_0^{1,p}(\Omega; W_p(\delta))$  such  that:
\begin{enumerate}\item $E (u)= u \quad $ a.e. in $\Omega_\eta$
\item There exists some  positive  number $C=C(\eta)$  such that  for  any 
$u\in W_0^{1,p} ( {\Omega_\eta}; W_p(\delta))$
$$ \| |\nabla E( u)| \|_{L^p(\Omega; W_p(\delta))}
 \le  C\Big(  \||\nabla u\||_{L^p (\Omega_{\eta/2}; W_p(\delta) )}+ \|u\|_{W^{1,p} (\Omega_\eta \setminus \overline{\Omega_{\eta/2}} ; W_p(\delta))}\Big).$$
\end{enumerate}
\end{lem}
Admitting this  for  the  moment, we prove  Theorem \ref{NCT2}.
First we  treat the  case  that $w\in P(\mathbf R_+)$. Then we  see that $W_p/F_\eta^p\in L^1(0,\eta)$.
In  fact we  have $\int_0^\eta W_p(t)/F_\eta(t)^p\,dt= \mu^{1-p}/(p-1).$
Now we assume  that
Hardy's inequality (\ref{2.11}) with $C=0$ does  not  hold. 
Then there exists a sequence of  functions $\{u_k\} \subset W_0^{1,p}(\Omega,W_p(\delta))\cap C(\Omega)$ satisfying
\begin{equation} \begin{cases}&\lim_{k\to \infty}\left(\int_{\Omega_\eta} |\nabla u_k|^pW_p(\delta)\,dx +   \int_{\Sigma_\eta}|u_k|^p W_p(\delta)\,d\sigma_\eta\right)=0,\\ &
 \int_{\Omega_\eta} {|u_k|^pW_p(\delta)}/{F_\eta^p(\delta)}\,dx =1 \quad (k=1,2,\cdots). 
 \end{cases} \label{6.1}
\end{equation}
By $w_k=u_k|_{\Omega_\eta}$ we denote the restriction of $u_k$ to $\Omega_\eta$.
Then $ E(w_k) \in W_0^{1,p}(\Omega;W_p(\delta))$ for $k=1,2,\ldots$. 
On  the  other hand, it  follows  from (\ref{6.1}) that
$w_k=u_k\to C$ a.e. in $\Omega_\eta$   for  some  constant  $C$ as $k\to\infty$. From (\ref{6.1}) we  have $C=0$. Then, by the  assumption 1 and the continuity of the  trace operator $T: W^{1,p}(\Omega\setminus \overline{\Omega_\eta}; W_p(\delta)) \mapsto 
 L^p(\Sigma_\eta; W_p(\delta ) )$  for  a  small  $\eta>0$,  we  have
 \begin{align} 1 &\le  \int_\Omega \frac{|E(w_k)|^pW_p(\delta)}{F_\eta^p(\delta)}\,dx \label{6.2}\\& \le  \gamma^{-1}\left(
 \int_\Omega |\nabla E(w_k)|^pW_p(\delta)\,dx  + L' \int_{\Sigma_\eta} |w_k|^pW_p(\delta)\,d\sigma_\eta  \right) \notag\\ 
& \le    \gamma^{-1}C'\Big(  \||\nabla u_k| \| _{L^p (\Omega_{\eta/2}; W_p(\delta) )}
+ \textcolor{black}{  \|u_k\|_{{W^{1,p}}(\Omega_\eta \setminus \overline{\Omega_{\eta/2}};W_p(\delta))} }\Big)^{\textcolor{black}{p}} \to 0, \quad\text{ as } k\to\infty,\notag
 \end{align}
 where $C'$ is  some  positive number  independent of  each $u_k$.
But  this  is a contradiction.

\par Secondly  we treat  the case that $w\in Q(\mathbf R_+)$.   Since
$\int_0^\eta W_p(t)/F_\eta(t)^p\,dt= +\infty$, we  see that $W_p/F_\eta^p\notin L^1(0,\eta)$. 
There is a sequence of  functions $\{u_k\} \subset W_0^{1,p}(\Omega; W_p(\delta))\cap C(\Omega)$ satisfying
\begin{align} &\lim_{k\to \infty}\int_{\Omega_\eta} |\nabla u_k|^pW_p(\delta)\,dx =0, \quad
 \int_{\Omega_\eta} \frac{|u_k|^pW_p(\delta)}{F_\eta^p(\delta)} \,dx =1 \quad (k=1,2,\cdots). \label{6.3}
\end{align}
As before  we  see that $u_k\to C$ (constant) in  $\Omega_\eta$. From (\ref{6.3}) we  have $C=0$.
Then, by the  assumption 1 and  the continuity of the  trace operator $T$, we  have
 \begin{align} 1 &\le  \int_\Omega \frac{|E(w_k)|^pW_p(\delta)}{F_\eta^p(\delta)}\,dx  \label{6.4}\\ &\le  \gamma^{-1}\left(
 \int_\Omega |\nabla E(w_k)|^pW_p(\delta) \,dx  - L' \int_{\Sigma_\eta} |w_k|^pW_p(\delta) \,d\sigma_\eta  \right)\\ 
& \le    \gamma^{-1}C'\Big( \| |\nabla u_k| \|_{L^p (\Omega_{\eta/2}; W_p(\delta) )}
+ \textcolor{black}{  \|u_k\|_{{W^{1,p}} (\Omega_\eta \setminus \overline{\Omega_{\eta/2}}; W_p(\delta))}} \Big)^{\textcolor{black}{p}}\to 0, \quad\text{ as } k\to\infty,\notag
\notag
 \end{align}
 where $C'$ is  some  positive number.
But  this  is a contradiction.
 \hfill $\Box$
\par\medskip
 \noindent{\bf Proof  of  Lemma \ref{l11}: } Since $\eta $  is  small and $\delta$ is  Lipschitz continuous, 
$\partial \Omega_\eta$  and $\partial \Omega_{\eta/2}$ are Lipschitz   compact manifolds. By  the  standard theory 
(  Theorem 1 in Section 5.4 of \cite{E}  for example ) we  have  an extension operator $\tilde{E} : W^{1,p} (\Omega_\eta \setminus \overline{\Omega_{\eta/2}} ; W_p(\delta))\mapsto
 W^{1,p}(\Omega\setminus \overline{\Omega_{\eta/2}}; W_p(\delta))$ such  that $\tilde E(u)= u$  a.e. in $\Omega_\eta \setminus \overline{\Omega_{\eta/2}}$, and 
$$ \| |\nabla \tilde{E}( u)|^p\|_{L^p(\Omega\setminus \overline{\Omega_{\eta/2}}; W_p(\delta))} 
\le  C(\eta) \textcolor{black}{ \|u\|_{{W^{1,p} }(\Omega_\eta \setminus \overline{\Omega_{\eta/2}};W_p(\delta))}}.$$
 Define  for  $u\in W_0^{1,p}(\Omega_\eta; W_p(\delta))$
  \begin{equation} { E}(u) = u, \, (x\in \overline{\Omega_{\eta/2}},);  \quad  \tilde{E}(u), \,  (x  \in  \Omega\setminus \overline{\Omega_{\eta/2}} ).
 \end{equation}
 Then the  assertion follows.  \hfill$\Box$
 
 \section{Proofs  of  Propositions \ref{ct1}  and  \ref{CT1}}
Proposition  \ref{ct1} is    known in a more general  fashion. In  fact  
a variant  is seen in Maz'ya   \cite{Ma} ( Lemma 2, p144). For  the sake  of  reader's convenience we  give an elementary verification. 
We note  that 
Proposition \ref{CT1} is  a direct  consequence  of  Proposition  \ref{ct1}.
\par\medskip
\noindent{\bf Proof of Proposition \ref{ct1}:} 
 For $w\in P(\R_+)$ and $\overline\varepsilon \in (0,\eta/2)$,  define 
\begin{equation} \varphi_{\overline\varepsilon}= \,\,
0 \,\,\,( 0\le t\le \overline\varepsilon) ;  \,\,\, \frac {f_\eta(\overline\varepsilon)-f_\eta(t)  }{f_\eta(\overline\varepsilon)-f_\eta(\eta/2) }\, \,\,(\overline\varepsilon \le t\le \eta/2);\,\,\,  1\, \,\,(\eta/2 \le t \le \eta). \label{7.1}\end{equation}
Noting  that $f_\eta(\overline\varepsilon)= \mu+ \int_{\overline\e}^\eta 1/w(s)\,ds\to +\infty$ as $\overline\varepsilon\to +0$, we  see  that
$ \int _0^\eta|\varphi_{\overline\varepsilon} '|^p W_p(t)\,dt=( f_\eta(\overline\varepsilon)-f_\eta(\eta/2))^{1-p} \to 0 \, \mbox{ as } \overline\varepsilon\to +0.$ 
On  the other hand  we  have  $\varphi_{\overline\varepsilon}(0)=0, \varphi_{\overline\varepsilon} (1)=1$ and  hence the  assertion is now  clear. 
Further we  note  that 
$$\int_{\overline\varepsilon}^\eta |\varphi_{\overline\varepsilon}|^p \frac{W_p(t)}{F_\eta(t)^p}\,dt
\ge  \int_{\eta/2}^\eta \frac{W_p(t)}{F_\eta(t)^p}\,dt=
\frac{f_\eta(\eta)^{1-p}- f_\eta(\eta/2)^{1-p}}{p-1}  >0
\quad  \mbox{ as } \overline\varepsilon\to +0.$$
\qed

\noindent{\bf  Proof  of  Proposition \ref{CT1}:}
If  a positive  number $\eta_0$  is sufficiently  small, then one  can  assume  that  $\delta\in C^2(\Omega_{\eta_0})$, $|\nabla \delta|=1$ in $ \Omega_{\eta_0}  $ and 
a manifolds $\{ x\in \Omega;  \delta =\eta\}$ is of  $C^2$  class for  $\eta \in (0,\eta_0]$. 
Let $\varphi_{\overline\varepsilon}$ be defined by (\ref{7.1}).
By  virtue  of  (\ref{4.2}) we  have
$$ \int_{\Omega_\eta}|\nabla \varphi_{\overline\varepsilon}(\delta(x))|^p {W_p(\delta(x))}\,dx\le 
 \int_{\Sigma} \,d\sigma \int_0^\eta \left | \varphi'_{\overline\varepsilon}(t) \right |^p {W_p(t)}
 (1+ct)  \,dt,$$
hence  the  assertion follows  from  Proposition \ref{ct1}.\qed

\section{Appendix}
\subsection{ Proof of Lemma \ref{lemma4.4}}
Let $v\in  G((0,\eta])$.
We  recall  that  
$$  g_\eta=( p' f_\eta)^{{1}/{p'}},\quad F_\eta= w f_\eta, \quad G_\eta= \mu+ \int_t^\eta\frac{ds}{F_\eta(s)}\quad\mbox{and}\quad  u= g_\eta \nu.$$
Let  us  prepare more  notations.
Let us  set 
\begin{equation}  v^{{p}/{2}}=z,\quad   z= a\varphi \quad\mbox{and}\quad a=(G_\eta)^{1/2},
\end{equation}
where  $ \varphi $ is some  function in $ C^1((0,\eta])$ with $\varphi(0)=0$.\par\noindent
Then we  have 
\begin{align*} z' &= \varphi' a + \varphi a'\quad\mbox{and}\quad 
|z'|^2 = (\varphi')^2 a^2 + \frac 12 (\varphi^2)' (a^2)' + \varphi^2 (a')^2.
\end{align*}
Noting  that $a^2= G_\eta,
 \,(a^2)' = -{1}/{F_\eta}\, ${and}  $\,  a' = - {1}/({2 F_\eta a}),$
we  have
\begin{align*} |z'|^2F_\eta= (\varphi')^2F_\eta G_\eta -\frac12 (\varphi^2)' + \frac{\varphi^2}{ 4F_\eta G_\eta} \ge -\frac12 (\varphi^2)' + \frac{\varphi^2}{ 4F_\eta G_\eta},
\end{align*}
and we  also have 
$$\varphi^2=\frac{v^p}{a^2}= \frac{ u^p}{g_\eta^pG_\eta} 
= \textcolor{black}{ \frac{1}{(p')^{p-1}}}\frac{u^p W_p}{ F_\eta^{p-1} G_\eta}. $$
As a result, we  have  the  following that is valid  for  any measurable  set $S\subset[0,\eta]$.
\begin{lem}\label{lemma8.1} Let $S$  be a measurable  set contained in $[0,\eta]$.
Assume that  $\eta>0$, $\mu>0$ and  $ w\in W(\mathbf R_+)$. For  any nonnegative $v\in  G([0,\eta])$
\begin{align} \int_S|(v(t)^{p/2})'|^2 F_\eta(t)\,dt
\ge -\frac12 \int_S(\varphi(t)^2)' \,dt+ \frac{1}{4(p')^{p-1}}\int_S\frac{u(t)^p W_p}{ F_\eta(t)^p G_\eta(t)^2}\,dt.
\end{align}
\end{lem}
\par\noindent{\bf Proof of Lemma \ref{lemma4.4}: } 
\textcolor{black}{
By using Lemma \ref{lemma8.1} with $S= [0,\eta]$ we immediately have (\ref{4.32'}).} Then we proceed to (\ref{4.32}).
We  assume  that $A=A(v,M)$ and $B=B(v,M)$ are defined by Definition \ref{df4.3}.
As for  the  first term of the right-hand side, we  have
\begin{align*}
-\frac12 \int_A(\varphi(t)^2)' \,dt &=-\frac12  \int_0^\eta(\varphi(t)^2)'\,dt +\frac12 \int_B(\varphi(t)^2)' \,dt\\&
=-\frac12 \varphi(\eta)^2 -\int_B\frac{z(t)^2 a'(t)}{a(t)^3}     +\int_B\frac{z(t) z'(t)}{a(t)^2}.
\end{align*}
Moreover we see that:
\begin{enumerate}
\item
 \begin{equation}\varphi(\eta)^2=\frac{z(\eta)^2}{a(\eta)^2}=\frac{1}{\mu}
\frac{|u(\eta)|^p}{g_\eta(\eta)^p} = \frac{|u(\eta)|^p}{\mu(p')^{p-1} f_\eta(\eta)^{p-1}}.
\end{equation}
\item
\begin{equation}-\int_B\frac{z(t)^2 a(t)'}{a(t)^3}\,dt=
\frac12 \int_B\frac{z(t)^2 }{a(t)^4 F_\eta(t)}\,dt=\frac{1}{2(p')^{p-1}}\int_B\frac{u(t)^p W_p(t)}{ F_\eta(t)^p G_\eta(t)^2}\,dt.
\end{equation}
\item
Noting  that $ zz'= p v^{p-1}v'/2 $, $ |v|\le {p'}|v'|F_\eta /{M}$ in $B$ and $G_\eta\ge \mu$, we  have
\begin{align*}\left|\int_B\frac{z(t)z'(t) }{a(t)^2 }\,dt\right|
&\le
\frac{p}{2}\int_B \frac{\left(\frac{p'}{M} |v'(t)| F_\eta(t)\right)^{p-1}|v'| }{G_\eta(t)} \,dt\\
& \le \frac{p}{2\mu} \left(\frac{p'}{M} \right)^{p-1}
\int_B  |v'(t)| ^{p}F_\eta(t)^{p-1}\,dt.
\end{align*}

\end{enumerate}
Finally  we get the desired inequality (\ref{4.32}) which proves  Lemma \ref{lemma4.4}.
\hfill $\Box$

\subsection{Auxiliary inequalities  in  the noncritical case}
If  we restrict ourselves to  the  case  that $w\in Q(\mathbf R_+)$, then Hardy's inequality (\ref{8.10})
follows from the next simple lemma which  is provided  in \cite{BM} for $p=2$ and $W_p=w=1$.
\begin{lem}\label{l2} \mbox{\rm ( $w\in Q(\mathbf R_+)$)}
Assume  that \textcolor{black}{ $ f\in C^{}([0,\eta])\cap C^1((0,\eta])$ } is a monotone nondecreasing function such  that $f(\eta)\le 1$.  Assume that $1<p<\infty$ and $w\in Q(\mathbf R_+)$.
Then for  every
  $u\in C_c^1((0,\eta])$, we  have 
\begin{equation} \int_0^\eta\left ( |u'(t)|^p -\Lambda_{p} \frac{|u(t)|^p}{F_\eta(t)^p} \right)W_p(t)\,dt \ge  
\int_0^\eta \left( |u'(t)|^p -\Lambda_{p} \frac{|u(t)|^p}{F_\eta(t)^p} \right)W_p (t)f (t) \,dt. \label{3.2}
\end{equation}
In particular we  have 
\begin{equation} \int_0^\eta |u'(t)|^p W_p(t)\,dt \ge  
\Lambda_{p} \int_0^\eta  \frac{|u(t)|^p}{F_\eta(t)^p} W_p (t) \,dt.  \label{8.10}
\end{equation}

\end{lem}
 \par\medskip
\noindent{\bf Proof of Lemma \ref{l2}:} Without loss of generality 
we  assume  that $f\ge 0$,  $f(\eta)=1$, and $u\ge0$.  Define  $ g= 1-f$. Then  $g\ge 0$  and  $ g'\le 0$.
Noting that $u\in C^1_c((0,\eta])$  and   $$ \frac d{dt}\left(\int_0^t\frac 1{w(s)}\,ds \right)^{1-p}= (1-p)\frac{W_p(t)}{F_\eta(t)^p},$$
by integration by parts we have 
\begin{align*}  & (p-1)\int_0^\eta  \frac{u(t)^p W_p(t)}{F_\eta(t)^p} g\,dt \\ &=
\int_0^\eta u(t)^p \left( \int_0^t  \frac1{w(s)}\,ds\right)^{1-p} g'(t)\,dt +  p \int_0^\eta u(t)^{p-1}  u'(t)  g(t)\left( \int_0^t  \frac1{w(s)}\,ds\right)^{1-p}\,dt.
\end{align*}
Since $g'=-f'\le 0$  and  $g\ge 0$,
\begin{align*} \frac1{p'} \int_0^\eta \frac{u(t)^p W_p(t)}{F_\eta^p(t)}  g(t)\,dt
\le \int_0^\eta u(t)^{p-1}  u' (t) g(t) \left( \int_0^t  \frac1{w(s)}\,ds\right)^{1-p}\,dt.
\end{align*}
By  H\"older's inequality, we have 
\begin{equation*}\begin{split}
\int_0^\eta u(t)^{p-1}  u' (t) g(t) &\left( \int_0^t  \frac1{w(s)}\,ds\right)^{1-p}\,dt\\
&\le  \left(  \int_0^\eta \frac{u(t)^p W_p(t)}{F_\eta(t)^p}  g(t)\,dt \right)^{1/{p'}}
 \left( \int_0^\eta  |u'(t) |^p
W_p (t)g(t)\,dt\right)^{1/p}.
\end{split}
\end{equation*}
Hence we  have 
\begin{align*} \frac1{p'}\left(  \int_0^\eta \frac{u(t)^p W_p(t)}{F_\eta(t)^p}  g(t)\,dt \right)^{1/p}
\le \left( \textcolor{black}{\int_0^\eta} |u'(t) |^p
W_p (t)g(t)\,dt\right)^{1/p}.
\end{align*}
Using $ g=1-f$ and  the definition of  $\Lambda_{p}$, we have (\ref{3.2}).
\hfill $\Box$
\par\bigskip
\noindent{{\bf Acknowledgments}}\par\medskip
The author would like to thank the referee for the precise advice in completing this paper.

\bigskip\bigskip
{\large\bf Department of Mathematics\\Faculty of Science \\ Ibaraki University\\
Mito, Ibaraki, 310, Japan}\\
\\
e-Mail; toshio.horiuchi.math@vc.ibaraki.ac.jp

\end{document}